\documentclass[11pt]{article}   
\usepackage{amssymb,amscd,latexsym}   
\usepackage{amsmath}
\usepackage{amsthm}
\usepackage{amsfonts}
\usepackage{mathdots}
\usepackage{mathbbol}
\usepackage[pagebackref]{hyperref}
\usepackage{color}
\usepackage{xcolor}
\usepackage{longtable}
\usepackage{lscape}
\usepackage{scalefnt}
\usepackage{setspace}
\usepackage{enumitem}
\usepackage{verbatim} 
\input xypic
\xyoption{all}

\newcommand{\rar}{\rightarrow}
\newcommand{\lar}{\longrightarrow}

\newcommand{\llar}{-\kern-5pt-\kern-5pt\longrightarrow}
\newcommand{\lllar}{-\kern-5pt-\kern-5pt\llar}
\newcommand{\surjects}{\twoheadrightarrow}
\newcommand{\kk}{\mathbb{k}}
\newcommand{\fm}{\mathfrak{m}}

\newtheorem{Theorem}{Theorem}[section]
\newtheorem{Lemma}[Theorem]{Lemma}
\newtheorem{Setup}[Theorem]{Setup}
\newtheorem{Corollary}[Theorem]{Corollary}
\newtheorem{Proposition}[Theorem]{Proposition}

\newtheorem{Conjecture}[Theorem]{Conjecture}

\newtheorem{Remark}[Theorem]{Remark}
\newtheorem{Example}[Theorem]{Example}
\newtheorem{Definition}[Theorem]{Definition}

\def\demo{\noindent{\bf Proof. }}
\def\depth{\mbox{\rm depth }}

\def\hht{{\rm ht}\,}
\def\Hom{\mbox{\rm Hom}}

\def\m{\mathfrak{m}}

\def\QED{\hfill$\Box$}
\def\qed{\QED}

\def\rank{\mbox{\rm rank}\,}
\def\rk{\mbox{\rm rank}\,}

\def\P{{\mathbb P}}

\def\YY{{\mathbf Y}}

\def\xx{{\mathbf x}}

\def\ff{{\mathbf f}}

\def\g2{{\mathbf g}}

\def\xx{{\mathbf x}}
\def\aa{{\mathbf a}}
\def\bb{{\mathbf b}}
\def\cc{{\mathbf c}}
\def\vv{{\mathbf v}}
\def\yy{{\mathbf y}}
\def\zz{{\mathbf z}}

\def\ee{{\bf e}}
\def\ss{{\bf s}}
\def\uu{{\bf u}}
\def\tt{{\bf t}}
\def\gg{{\bf g}}
\def\ddd{{\bf d}}

\def\m{{\mathfrak m}}

\begin{document}

\begin{center}
	
	{\Large{\bf\sc Equigenerated  Gorenstein ideals of codimension three}}
	\footnotetext{AMS 2010 Mathematics
		Subject Classification (2010   Revision). Primary  13H10, 13D02; Secondary 13A30, 13C13, 13C15, 13C40.}

\bigskip

{\sc Dayane Lira}\footnote{Under a PhD fellowship from CAPES, Brazil (88882.440720/2019-01)}
\quad\quad
{\sc Zaqueu Ramos}\footnote{Under a post-doc fellowship from INCTMAT/Brazil (88887.373066/2019-00)}
\quad\quad
{\sc Aron  Simis}\footnote{Partially supported by a CNPq grant (302298/2014-2)}

\end{center}

	
\begin{center}
\end{center}



\begin{abstract}
\noindent

We focus on the structure of a homogeneous Gorenstein ideal $I$ of codimension three in a standard polynomial ring $R=\kk[x_1,\ldots,x_n]$ over a field $\kk$, assuming that $I$ is generated in a fixed degree $d$. For such an ideal $I$ this degree comes along with the minimal number of generators of $I$ and the degree of the entries of the associated skew-symmetric matrix in a simple formula. 
We give an elementary characteristic-free argument to the effect  that, for any such data linked by this formula, there exists a Gorenstein ideal $I$ of codimension three filling them.
We conjecture that, for arbitrary $n\geq 2$, an ideal $I\subset \kk[x_1,\ldots,x_n]$ generated by a general set of $r\geq n+2$  forms of degree $d\geq 2$ is Gorenstein if and only if $d=2$ and $r= {{n+1}\choose 2}-1$.
We prove the `only if'  implication of this conjecture when $n=3$.  For arbitrary $n\geq 2$, we prove that if $d=2$ and $r\geq (n+2)(n+1)/6$ then the ideal is Gorenstein if and only if $r={{n+1}\choose 2}-1$, which settles the `if' assertion of the conjecture for $n\leq 5$. Finally, we elaborate around one of the questions of Fr\"oberg--Lundqvist.
In a different direction, we reveal a connection between the Macaulay inverse and the so-called Newton dual, a matter so far not brought out to our knowledge.
Finally, we consider the question as to when the link $(\ell_1^m,\ldots,\ell_n^m):\mathfrak{f}$ is equigenerated, where $\ell_1,\ldots,\ell_n$ are independent linear forms and $\mathfrak{f}$ is a form, is given a solution in some important cases.

\end{abstract}


\section*{Introduction}

The literature on Gorenstein ideals is vast (for a tiny list, see \cite{BE1}, \cite{Wat}, \cite{Diesel}, \cite{Harima}, \cite {Iarr-Ems}, \cite{Iarr}, \cite{Iarrobino_book}, \cite{KU}, \cite{Rag-Zap}, \cite{Ia-Sri}, \cite{Kleppe}, \cite{Kustin}, \cite{EliasRossi}, \cite{trio_Gor}, \cite{Jeli}, \cite{Roig}), but perhaps not so much with an emphasis on the equigenerated case, even when the codimension is three.
Largely, one expects  this case to be easier to handle, by doing away with some of the technical difficulties of the arbitrary homogeneous case, where the  numerical data are more involved.
The question, of course, is as to whether there is a net gain in this restriction. We hope that the overall simplicity of the arguments in this work will justify doing it.


By and large we have been initially inspired by the classification ideas of \cite{Kustin} and \cite{abc}.
A fairly understood case is that of an ideal $I\subset R=\kk[x,y,z]$ of finite colength generated by quadrics. 
For example, one has:

\smallskip

\noindent {\bf Theorem.}{\rm (\cite[Theorem 2.1 (i) and Theorem 3.6]{abc})}
Let $\kk$ be a field of characteristic $\neq 2$ and let $I\subset R=\kk[x,y,z]$ be an ideal of finite colength minimally generated by five quadrics. Then the following conditions are equivalent:

	(i) $I$ is a Gorenstein  ideal.
	
    (ii)  $I$ is syzygetic.
    
	(iii) There exist $3$ independent linear forms $\{l_1,l_2,l_3\}\subset \kk[x,y,z]$ such that 
	$$I=(xl_2,xl_3,yl_3, xl_1-yl_2,xl_1-zl_3).$$	

Thus, in the case of quadric generators in dimension three, it gives characterizations other than the  Pfaffian recipe of Buchsbaum--Eisenbud.
Condition (ii) and (iii) may suggest that  such Gorenstein ideals live on a dense Zariski open set in the parameter space of the coefficients of a set of generators.
In fact, it turns out as we prove more generally, that for arbitrary $n\geq 2$, an ideal $I\subset \kk[x_1,\ldots,x_n]$ generated by a general set of $r\geq (n+2)(n+1)/6$  quadrics  is Gorenstein if and only if $r={{n+1}\choose 2}-1$.

As a matter of fact, we conjecture an encompassing result to the effect  that, for arbitrary $n\geq 2$, an ideal $I\subset \kk[x_1,\ldots,x_n]$ generated by a general set of  $r\geq n+2$ forms of degree $d\geq 2$ is Gorenstein if and only if $d=2$ and $r= {{n+1}\choose 2}-1$. 
The previously mentioned theorem proves one implication of this conjecture. We are able to settle the reverse implication in dimension three, for which we give two different proofs.
The first is partially based on previously existing material, but we state all facts {\em ab initio} in order to have a coherent expos\'e. Thus, one will recognize in the background some of the results by Hochster--Laksov (\cite{HoLa}), Diesel (\cite{Diesel}), Anick (\cite{Anick}) and Migliore--Mir\'o-Roig (\cite{Migl-Roig}).
The second proof  is achieved by making explicit that all equigenerated Gorenstein ideals of codimension $3$, of given degree $d\geq 3$ and $r\geq 5$ number of generators, are parameterized by a proper closed subset of the space of parameters. For this, we bring up 
 estimating the socle degree, for which it became necessary to give the details of the graded free resolution. Due to the length and technicality of the involved arguments, we deferred the details to the Appendix.
 
An aftermath is a discussion around \cite[Question 2.5]{Fro-Lund} that asks for the nature of the coefficients of the power series $((1-t^d)^r/(1-t)^n)_+$ for $r\geq n$.

In another direction, the minimal number $\mu(I)$ of generators of a codimension $3$ Gorenstein ideal $I\subset \kk[x,y,z]$ generated in degree $d$ obeys a certain additional constraint, besides being an odd integer. 
This obstruction is formalized in terms of the available  numerical data, to wit, a pair $(d,r)$ of integers such that $d\geq 2$ and $r\geq 3$ will be called an {\em equigenerated Gorenstein  virtual datum} in dimension $n\geq 3$ if 
$(r-1)/2$ is a factor of $d$.
Such a datum will be said to be {\em proper} (in dimension $n$) if there exists a codimension $3$ Gorenstein ideal in $\kk[x_1,\ldots,x_n]$ generated by $r$ forms of degree $d$, and moreover, the entries of a corresponding $r\times r$ skew-symmetric matrix $\phi$  generate an ideal of codimension $n$. Here,  $2d/(r-1)$ will be the degree of any entry of $\phi$. 
We prove that any virtual datum is proper (Theorem~\ref{properness}).
This result is characteristic-free. 
This question has been tackled before in \cite{Diesel} and by Conca--Valla (\cite{Con-Valla}) via an argument based upon a particular skew--symmetric matrix. An approach of the same sort, drawing upon the early example of Buchsbaum and Eisenbud (\cite[Proposition 6.2]{BE1}), will appear in the PhD thesis of the first author.

Another subject dealt with in this paper is what we call the $(x_1^m,\ldots, x_n^m)$-{\em colon problem}, for lack of better terminology.
It has long been known (see \cite[Proposition 1.3]{BE0}) that, in arbitrary characteristic, any homogeneous Gorenstein ideal of codimension $n$ in $\kk[x_1,\ldots,x_n]$ can be obtained as a colon ideal $(x_1^m,\ldots, x_n^m):\mathfrak{f}$, for some integer $m \geq 1$ and a form $\mathfrak{f}$. 
In this regard, two questions naturally come up: first, is there a more definite relation between $\mathfrak{f}$ and the ideal $I$? Second, is there a characterization as to when the resulting Gorenstein ideal is equigenerated in terms of the exponent $m$ and the form $\mathfrak{f}$?

We give an answer to the first question, at least in characteristic zero, in terms of the Macaulay inverse setup and the notion of Newton duality as introduced in \cite{CostaSi} and \cite{DoriaSi}.
This is the content of Proposition~\ref{NewtondualxMinverse} which says that Macaulay inverse to $I$ is the (socle-like) Newton dual of the form $\mathfrak{f}$. The proposition also gives a complete characterization of $\mathfrak{f}$, and its degree, under the condition that none of its nonzero terms belongs to the ideal $(x_1^m,\ldots, x_n^m)$.
The use of the Newton duality concept in these matters seem to be new, as far as we can tell.

The second problem above is more delicate.  One usefulness of solving it is yet another simple way of producing equigenerated Gorenstein ideals. Obviously, merely controlling degrees in both terms of the colon operation does not lead too far -- e. g., in characteristic zero the ideal $(x^d,y^d,z^d): (x+y+z)^{d-1}$,  for $d\geq 2$ , is $d$-equigenerated, while $(x^4,y^4,z^4): x^{3}+y^{3}+z^{3}$  is not equigenerated.
Yet, we solve this question in the case where $I$ has linear resolution (Theorem~\ref{colon_problem_linear_case}), in terms of a degree constraint and show that, when the directrix $\mathfrak{f}$ is a power of $x_1+\cdots +x_n$, the solutions lie on a dense open set of the space  of parameters (Proposition~\ref{sum_power_case}).
The main part of the latter requires  characteristic zero because it invokes a result of  Stanley, with algebraic proofs by Reid--Roberts--Roitman (\cite{RRR91} and Oesterl\'e in (\cite[Apendix A]{BuChSi})).

In the way of considering the last problem, we introduce the notion of {\em pure power gap}, an integer measuring how far off from the socle degree is an exponent of the powers of independent linear forms lying in the ideal.
We give the basic role of this invariant, hoping it will be useful in other contexts. 

Finally,  we look at some behavior of the usual algebras attached to an ideal, such as the Rees algebra and the fiber cone algebra. The relation between the depth of these algebras and the nature of the rational map defined by the linear space $I_d$ is discussed in Section~\ref{Rees_all} by completely elementary ways. 

A list of the main results: Theorem~\ref{properness},  Theorem~\ref{Old_Theorem 2.7}, Theorem~\ref{socle_always_in_degree2}, Theorem~\ref{colon_problem_linear_case}  and  Theorem~\ref{main_syzygetic}.
Other significant results:
Proposition~\ref{Soc1bis},
 Proposition~\ref{NewtondualxMinverse}, Proposition~\ref{sum_power_case}.and Proposition~\ref{power_gap_characterization}.,
Proposition~\ref{examplemonomialbis} and Proposition~\ref{examplemonomialbisbis}.  

\smallskip

{\sc Acknowledgment:} Upon posting  a first version of this work on the arXiv, it has been brought to our knowledge by A. Iarrobino that some of our results have been considered before.
We thank him   for pointing out the missing references, which are now included.
Since our approach has often different features from the previously existent literature, we decided to keep it as a possible new angle of consideration.
We also thank R. Fr\"oberg for useful conversations around his celebrated conjecture and correlated issues, where he pointed out some inaccuracies.


\section{Numerical data}\label{virtual_data_all}

\subsection{The socle}\label{socle_info}

Let $R=\kk[x_1,\ldots,x_n]$ be a polynomial ring over a field $\kk.$ Denote by $\fm$ the  maximal homogeneous ideal of $R.$ Given a homogeneous $\fm$-primary ideal $I\subset R$, consider the least integer $s$  such that $\fm^{s+1}\subset I.$ Then, the graded $\kk$-algebra $R/I$ can be written as
$$R/I=\kk\oplus (R/I)_1\oplus(R/I)_2\oplus \cdots\oplus(R/I)_s$$  
with $(R/I)_s\neq 0.$ The {\em socle}  of $R/I,$ denoted ${\rm Soc}(R/I),$ is the ideal $I:\fm/I\subset R/I$. Since $I$ and $\m$ are homogeneous ideals, then  ${\rm Soc}(R/I)$ is a homogeneous ideal of $R/I.$ In particular,

$${\rm Soc}(R/I)={\rm Soc}(R/I)_1\oplus\cdots \oplus{\rm Soc}(R/I)_s.$$
The integer $s$ is the {\it socle degree} of $R/I.$
It is well known that $R/I$ is a Gorenstein ring if and only if ${\rm Soc}(R/I)={\rm Soc}(R/I)_s$ and $\dim_{\kk}{\rm Soc}(R/I)_s=1.$

Recall that a codimension $3$ Gorenstein ideal  $I\subset R=\kk[x_1,\ldots,x_n]$ is generated by the Pfaffians of an $r\times r$ skew-symmetric matrix $\mathfrak{S}$, with $r=\mu(I)$ odd.
When $I$ is equigenerated, say, in degree $d\geq 1$, then the columns of $\mathfrak{S}$ must be homogeneous of some standard degree $d_i,\, 1\leq i\leq r$.
By the nature of each generator of $I$ as a maximal Pfaffian of $\mathfrak{S}$, it immediately follows that
$$2d=d_1+d_2+\cdots +d_{r-2}+d_{r-1}=d_1+d_3+\cdots +d_{r-1}+d_{r} =\cdots = 
d_2+d_3+\cdots + d_{r-1}+d_{r}.$$
By an elementary argument,  $d_1=d_2=\cdots = d_r$.

It follows that $2d=(r-1)d'$, where $d'$ is the common values of the $d_i$'s, i.e.,
\begin{equation}\label{Gor_obstruction}
	d= \frac{r-1}{2}\, d'.
\end{equation}
Thus, an equigenerated Gorenstein ideal $I$ of codimension $3$ complies with this relation, where $r$ is the number of generators, $d$ their common degree and $d'$ is the degree of any entry of the skew--symmetric matrix defining $I$.

\begin{Lemma}\label{socle_degree_general}
	Let $I\subset R=\kk[x,y,z]$ denote an equigenerated codimension $3$ Gorenstein ideal with virtual datum $(d,(2d+d')/d')$, where $d'\geq 1$ is the skew degree of $I$ as introduced earlier.
	Then the socle degree of $R/I$ is $2d+d'-3$.
\end{Lemma}
\demo It is well-known that, pretty generally, the socle degree of a graded Artinian Gorenstein quotient $R/I$ of a standard graded polynomial ring $R=\kk[x_1,\ldots,x_n]$ is given by $D-n$, where $D$ is the last shift in the minimal graded $R$-resolution of $R/I$ (see, e.g., \cite[Lemma 1.3]{KU}).
In the present case, the resolution has length $3$ and the first syzygies have shifted degree $d+d'$.
Therefore, $D=d+d+d'=2d+d'$, hence the socle degree is as stated.
\qed

\smallskip

For reference convenience, we introduce the following terminology:
\begin{Definition}\label{virtual_data}\rm
	A {\em codimension $3$ equigenerated Gorenstein  virtual datum} in dimension $n\geq 3$ is a pair $(d,r)$ of integers such that:
	
	(i) $d\geq 2$ and $r\geq 3$.
	
	(ii) $(r-1)/2$ is a factor of $d$.
\end{Definition} 
If no misunderstanding arises, we will mostly omit `codimension $3$ equigenerated Gorenstein'	in the above terminology.
	Given such a datum, the integer $d':=2d/(r-1)$ will be called  the {\em skew-degree} of the datum.
	
	The datum is said to be {\em proper} (in dimension $n$) if there exists a codimension $3$ Gorenstein ideal $I$ in $\kk[x_1,\ldots,x_n]$ with this datum, satisfying $\hht I_1(\phi)=n$, where $\phi$ is the skew-symmetric matrix whose maximal Pfaffians generate $I$. In this case the skew-degree is in fact the degree of any entry in the skew-symmetric matrix $\phi$, perhaps justifying its designation.

Note that, since $n\geq 3$ then, rightfully,  $r\leq 2d+1$.
An important case is when $d'=1$, i.e., $r=2d+1$. It is colloquially  referred to as the {\em linear case}.
This happens, e.g., when $d$ is a prime number and $r\geq 4$.

\subsection{Retrieval from numerical data} 

In this part we prove that any virtual datum in dimension $n\ge 3$ is proper. The argument works in arbitrary characteristic. 
A very similar question has been dealt with  in \cite{Diesel}.
For $n=3$ and homogeneous Gorenstein ideals which are not necessarily equigenerated, this problem has been solved in \cite[Section 1]{Con-Valla}.
Our approach for the equigenerated case  and $n\geq 3$ is based on a  simple reparametrization-like device as follows.

Quite generally, for  a standard graded polynomial ring $\kk[z_1,\ldots,z_N]$ over an infinite field $\kk$, and 
an integer $p\geq 1$, consider the following injective  $\kk$-algebra map
$$\zeta_{p}:\kk[z_1,\ldots,z_N]\to \kk[z_1,\ldots,z_N], \quad z_i\mapsto z_{i}^{p}\,(1\leq i\leq N).$$

\begin{Lemma}\label{zetacomuta}
	Let $>$ be a monomial order on $\kk[z_1,\ldots,z_{N}].$ Then, ${\rm in}_{>}(\zeta_{p}(f))=\zeta_{p}({\rm in}_{>}(f))$
	for every form $f$.
\end{Lemma}
\demo Write $f=a_0\zz^{\vv_0}+\cdots+a_s\zz^{\vv_s}$,
where  $a_i\neq 0$ for every $i$ and $\zz^{\vv_0}>\zz^{\vv_j}$ for every $j\geq 0.$ In particular, ${\rm in}_{>}(f)=\zz^{\vv_0}.$ We have
$$\zeta_{p}(f)=a_0(\zz^{\vv_0})^p+\cdots+a_s(\zz^{\vv_s})^p.$$
By the properties of a monomial order, one has $(\zz^{\vv_0})^p>(\zz^{\vv_j})^p$ for every $j\geq 0.$ Thus, ${\rm in}_{>}(\zeta_{p}(f))=(\zz^{\vv_0})^p=\zeta_{p}(\zz^{\vv_0})=\zeta_{p}({\rm in}_{>}(f)).$
\qed

\smallskip

If $M$ is any matrix with entries in $\kk[z_1,\ldots,z_N]$, we denote by $\zeta_p(M)$ the matrix obtained by evaluating $\zeta_p$ at every entry of $M$.
We focus on the case where $r$ is a given odd integer and $X$ denotes  the $r\times r$ generic skew-symmetric matrix.
Let $B$ stand for the polynomial ring over $\kk$ in the nonzero entries of $X$.
\begin{Lemma}\label{modelogenerico}
	For any integer $p\geq 1$ the $(r-1)$-Pfaffians of the skew-symmetric matrix  $\zeta_p(X)$ generate an ideal of $B$ of height $3$.
\end{Lemma}
\demo It is well-known that  an $(r-1)$-Pfaffian is a polynomial in the entries of the source skew-symmetric matrix (see, e.g., \cite[Proposition 159]{Muir}).
 Since $\zeta_p$ is a homomorphism of $\kk$-algebras, then the ideal   generated by the $(r-1)$-Pfaffians of $\zeta_p(X)$ is $\zeta_p(I)$, where $I$ is the ideal of $(r-1)$-Pfaffians of $X$.  Thus, in order to conclude the proof it suffices to show that $\hht \zeta_p(I)\geq 3.$ 

By \cite[Theorem 5.1]{HT} there is a monomial order $>$ on $B$ such that the $(r-1)$-Pfaffians of $X$ constitute a Gr\"obner base of $I.$ In particular, if $f_1,\ldots,f_r$ are the $(r-1)$-Pfaffians of $I$ then ${\rm in}_{>}(I)=({\rm in}_{>}(f_1),\ldots,{\rm in}_{>}(f_r)).$ Thus,
\begin{eqnarray}\label{inzeta}\nonumber
	{\rm in}_{>}(\zeta_p(I))&\supset& ({\rm in}_{>}(\zeta_p(f_1)),\ldots,{\rm in}_{>}(\zeta_p(f_r)))\\ 
	&=&(\zeta_p({\rm in}_{>}(f_1)),\ldots,\zeta_p({\rm in}_{>}(f_r)))\quad(\mbox{by Lemma~\ref{zetacomuta}})\\ 
	&=&({\rm in}_{>}(f_1)^p,\ldots,{\rm in}_{>}(f_r)^p)\nonumber
\end{eqnarray}
and, consequently,
\begin{eqnarray*}
\quad\quad\quad\quad\quad  	\hht\zeta_p(I)&=&\hht {\rm in}_{>}(\zeta_p(I))\\
	&\geq& \hht ({\rm in}_{>}(f_1)^p,\ldots,{\rm in}_{>}(f_r)^p)\quad(\mbox{by \eqref{inzeta}})\\
	&=& \hht ({\rm in}_{>}(f_1),\ldots,{\rm in}_{>}(f_r)) 
	=\hht {\rm in}_{>}(I)\\
	&=&\hht I=3.\quad \quad\quad\quad\quad\quad \quad \quad\quad\quad\quad\quad\quad \quad\quad\quad\quad\quad   \quad\quad\quad\quad\quad  \quad       \square
\end{eqnarray*}

\begin{Theorem}\label{properness} Let $n\geq 3$ be an integer. Then every virtual datum $(d,r)$ in dimension $n$ is proper.
\end{Theorem}
\demo Denote by $d'$ the skew degree of the virtual datum $(d,r).$ Define $u={r\choose 2 }-n\leq {r\choose 2 }-3.$ We induct on $u.$ If $u=0$, by  Lemma~\ref{modelogenerico} the ideal generated by the $(r-1)$-Pfaffians of the matrix $\zeta_{d'}(X)$
is a Gorenstein ideal with datum $(d,r).$ 

Now, suppose that the result is true for a certain $1\leq u<{r\choose 2 }-3.$
Then there is an $r\times r$ skew-symmetric matrix $\Psi=(g_{i,j})$ whose entries $g_{i,j}$ are forms of degree $d'$  in $\kk[x_1,\ldots,x_n]$, such that $\hht I_1(\Psi)=n$ and $I={\rm Pf}_{r-1}(\Psi)$ is a codimension 3 Gorenstein ideal. Since one can assume that $r$ is fixed, we may argue by descending induction on $n$ instead. Thus, we are assuming that $n>3$. Since $I$ is an unmixed homogeneous ideal of codimension $3$, then  by a homogeneous version of the prime avoidance lemma,  there is a linear form $\ell\in R= \kk[x_1,\ldots,x_n]$ that is $R/I$-regular. Without loss of generality, we can suppose that $\ell=x_n-\sum_{i=1}^{n-1}\alpha_ix_i.$ Now consider the following surjective $\kk$-algebra homomorphism:
$$\pi:\kk[x_1,\ldots,x_n]\surjects \kk[x_1,\ldots,x_{n-1}],\quad x_{i}\mapsto x_i\,(1\leq i\leq n-1),\,x_{n}\mapsto \sum_{i=1}^{n-1}\alpha_ix_i. $$
Using again the fact that Pfaffians are polynomials in the entries of the source matrix, $\pi$ induces a $\kk$-algebra isomorphism
$$\kk[x_1,\ldots,x_n]/({\rm Pf}_{r-1}(\Psi),\ell)\simeq \kk[x_1,\ldots,x_{n-1}]/{\rm Pf}_{r-1}(\widetilde{\Psi}) $$
and
$$\kk[x_1,\ldots,x_n]/(I_1(\Psi),\ell)\simeq \kk[x_1,\ldots,x_{n-1}]/I_1(\widetilde{\Psi})$$
where $\widetilde{\Psi}=(\pi(g_{ij})).$
By the second isomorphism above,  $I_1(\widetilde{\Psi})$ has height one. Since $\ell$ is  regular on $\kk[x_1,\ldots,x_n]/I$ then $\kk[x_1,\ldots,x_{n-1}]/{\rm Pf}_{r-1}(\widetilde{\Psi})$ is a Gorenstein ring of codimension $$(n-1)-\dim \kk[x_1,\ldots,x_{n-1}]/{\rm Pf}_{r-1}(\widetilde{\Psi})=(n-1)-(n-4)=3.$$
Thus, ${\rm Pf}_{r-1}(\widetilde{\Psi})$ is a codimension $3$ Gorenstein ideal in $\kk[x_1,\ldots,x_{n-1}]$ with datum $(d,r)$, satisfying the condition $\hht I_1(\widetilde{\Psi})=1$.
\qed

\section{Parametrization}

\subsection{The span problem}\label{Parameters}

Let $R=\kk[\xx]=\kk[x_1,\ldots,x_n]$ denote a polynomial ring over a field $\kk$.
Fix an integer $d\geq 1$ and let $\ff:=\{f_1,\ldots,f_r\}\subset R$ be a set of forms of degree $d$.

Set $f_{t}=\sum_{|\alpha|=d} \lambda_{\alpha}^{(t)}\xx^{\alpha}, 1\leq t\leq r$, and $N=\dim_{\kk} R_{d}-1={d+n-1\choose d}-1$ and consider the  parameter map $R_d\times\cdots\times R_d \to (\P^{N})^r$ that associates to the vector $[f_1\,\cdots\,f_r]$ the point
$$P_{\ff}:=(\lambda_{d,0,\ldots,0}^{(1)}:\cdots:\lambda_{0,\ldots,0,d}^{(1)};\ldots;\lambda_{d,0,\ldots,0}^{(r)}:\cdots:\lambda_{0,\ldots,0,d}^{(r)}).$$
We will harmlessly use the same notation $\ff$ for either the set or the vector of the given forms.

Fix an integer $e\geq 0$ and let $D:=\dim_{\kk}(R_{d+e})$. The   $\kk$-vector subspace $R_{e}f_1+\cdots+R_{e}f_{r}$ of $R_{d+e}$ is spanned by the set 
$\{\xx^{\alpha}f_1\,|\,|\alpha|=e\}\cup\cdots\cup\{\xx^{\alpha}f_{r}\,|\,|\alpha|=e\}$.
The  coefficient matrix of this set with respect to canonical basis of $R_{d+e}$ is a $D\times (r\dim_{\kk} R_{e})$ that can be written in the following form:
\begin{equation}\label{parameter_matrix}
	\left(\begin{array}{c|c|ccccccccc}
M_{1}&\ldots&M_{r}
\end{array}
\right),
\end{equation}
where $M_{t}$ is the coefficient matrix of the set $\{\xx^{\alpha}f_t\,|\,|\alpha|=e\}$ with respect to the canonical basis of $R_{d+e}.$  In particular, the entries of the block $M_{t}$ involve only the coefficients  $(\lambda_{d,0,\ldots,0}^{(t)},\ldots,\lambda_{0,\ldots,0,d}^{(t)})$ in the $t$th copy of $\P^N$. 

We can rewrite these data in terms of  projective coordinates. For this, let $\P^N={\rm Proj}(\kk[Y_{d,0,\ldots,0}^{(t)},\ldots,Y_{0,\ldots,0,d}^{(t)}])$ stand for  the $t$th copy of $\P^N$.
Let $M_{d,r,e}$ stand for the matrix obtained from (\ref{parameter_matrix})  upon replacing any entry 
$\lambda_{\alpha}^t$ in the block $M_t$  by the corresponding variable $Y_{\alpha}^t$.
Clearly, evaluating $M_{d,r,e}$ on the coordinates of the parameter point $P_{\ff}$ gives back the matrix in (\ref{parameter_matrix}), henceforth denoted $M_{d,r,e}(P_{\ff})$.

The following simple example may help visualizing the shape of  $M_{d,r,e}$.
\begin{Example}\rm Let $r=d=n=2$ and $e=1$.
	Then $D=4$ and one gets:  
	$$M_t=\left(\begin{array}{cccccccc}
	\medskip
	\lambda^{(t)}_{2,0}&0\\
	\medskip
	\lambda^{(t)}_{1,1}&\lambda^{(t)}_{2,0}\\
	\medskip
	\lambda^{(t)}_{0,2}&\lambda^{(t)}_{1,1}\\
	\medskip
	0&\lambda^{(t)}_{0,2}\\
	\end{array}\right) \:(t=1,2)
	\quad\mbox{and}\quad
	M_{2,2,1}=\left(\begin{array}{cc|cc}
	Y^{(1)}_{2,0}&0&Y^{(2)}_{2,0}&0\\ [2pt]
	Y^{(1)}_{1,1}&Y^{(1)}_{2,0}&Y^{(2)}_{1,1}&Y^{(2)}_{2,0}\\[2pt]
	Y^{(1)}_{0,2}&Y^{(1)}_{1,1}&Y^{(2)}_{0,2}&Y^{(2)}_{1,1}\\ [2pt]
	0&Y^{(1)}_{0,2}&0&Y^{(2)}_{0,2}\\
	\end{array}\right).
	$$ 
	\end{Example} 
In this example the matrix has maximal rank, but this will not be the case in arbitrary cases.



\smallskip

If $R_{e}f_1+\cdots+R_{e}f_{r}=R_{d+e}$ we loosely say that  the $d$-forms  $f_1,\ldots,f_r$ {\em span in degree} $d+e$.
A central question here asks when the inclusion $V(I_{D}(M_{d,r,e}))\subset (\P^{N})^r$ is proper, i.e., if the matrix $M_{d,r,e}$ has maximal rank. 
This is tantamount to having that, for given $r,d,e$, `most' tuples $(f_1,\ldots,f_r)$ span.

As a matter of further notation, for given data $\{d,r,e\}$ as above denote by $D_{d,r,e}$ the rank of the matrix $M_{d,r,e}.$
Accordingly, $V_{D_{d,r,e}}(M_{d,r,e})\subset (\mathbb{P}^{N})^r$ will denote the subvariety defined by the minors of $M_{d,r,e}$ of order  $D_{d,r,e}$. The following properties are immediately verified:

\begin{enumerate}[label=(\Alph*)]
	\item \label{A} $D_{d,r,e}\leq \min\{\dim_{\kk}R_{d+e},r\dim_{\kk}R_e\}.$ Equivalently, 
	\begin{equation}
	\dim_{\kk} R_{d+e}- D_{d,r,e}\geq \max\{\dim_{\kk} R_{d+e}-r\dim_{\kk} R_{e},\,0\}.
	\end{equation}
	\item \label{B} For every sequence  $\ff\in (R_d)^r,$ $\dim_{\kk} [(\ff)]_{d+e}=\rank M_{d,r,e}(P_{\ff})\leq D_{d,r,e}.$ Equivalently, 
	\begin{equation}
	\dim _{\kk}[R/(\ff)]_{d+e}\geq \dim_{\kk} R_{d+e}- D_{d,r,e}.
	\end{equation}
	\item  \label{C} For every sequence $\ff\in (R_d)^r$  such that $P_{\ff}\in(\mathbb{P}^N)^r\setminus V(I_{D_{d,r,e}}(M_{d,r,e})),$ 
	\begin{equation}
	\dim _{\kk}[R/(\ff)]_{d+e}= \dim_{\kk} R_{d+e}- D_{d,r,e}.
	\end{equation}
\end{enumerate}   

\begin{Definition}\rm
	Let $R=\kk[x_1,\ldots,x_n]$.
Given $d ,r,e$, a set of forms $\ff\in (R_d)^r$ is {\em $(d,r,e)$-extremal} if $$\dim_{\kk}[R/I]_{d+e}=\max\{\dim_{\kk} R_{d+e}-r\dim_{\kk} R_{e},\,0\}$$
 where $I=(\ff)\subset R$.
\end{Definition}

For convenience, we isolate some basic facts in the following lemma.

\begin{Lemma}\label{extremal_properties}
		Let $R=\kk[x_1,\ldots,x_n]$. Given integer data $d, r,e$ as above, one has:
		\begin{enumerate}
			\item[{\rm (i)}] If there is a $(d,r,e)$-extremal set in $(R_d)^r$ then  every set $\ff\in (R_d)^r$  such that $P_{\ff}\in(\mathbb{P}^N)^r\setminus V(I_{D_{d,r,e}}(M_{d,r,e}))$ is $(d,r,e)$-extremal.
			In particular, for any such $\ff$ one has
			$$\dim_{\kk} R_{d+e}- D_{d,r,e}= \max\{\dim_{\kk} R_{d+e}-r\dim_{\kk} R_{e},\,0\}.$$
			\item[{\rm (ii)}] For $e=0$, every $\ff\in (R_d)^r$ such that $P_{\ff}\in(\mathbb{P}^N)^r\setminus V(I_{D_{d,r,0}}(M_{d,r,0}))$ is $(d,r,0)$-extremal.
			\item[{\rm (iii)}] For $e=1$, every set $\ff\in (R_d)^r$ such that $P_{\ff}\in(\mathbb{P}^N)^r\setminus V(I_{D_{d,r,1}}(M_{d,r,1}))$ is $(d,r,1)$-extremal.
			\item[{\rm (iv)}] If $e>n(d-1)$ and $r\geq n$ then every $\ff\in (R_d)^r$ such that $P_{\ff}\in(\mathbb{P}^N)^r\setminus V(I_{D_{d,r,e}}(M_{d,r,e}))$ is $(d,r,e)$-extremal.
		\end{enumerate}
\end{Lemma}
\demo (i) This follows from properties \ref{A}, \ref{B} and \ref{C} above.

(ii) This is just linear independence over $\kk$.

(iii) This is the content of \cite[Theorem 1]{HoLa}.

(iv)  In fact, for every $e> n(d-1),$ $\max\{\dim_{\kk} R_{d+e}-r\dim_{\kk} R_{e},\,0\}=0.$ For a sequence $\ff_0\in( R_{d})^r$ such that $\{x_1^d,\ldots,x_n^d\}\subset \ff_0$ we have  $\dim_{\kk}[R/(\ff_0)]_{d+e}=0,$ that is, $\ff_0$ is $(d,r,e)$-extremal. Thus,  the claim follows by (i).
\qed

\medskip

Recall that a set $\ff=\{f_1,\ldots,f_r\}$ of forms of the same degree is {\em general} when the corresponding parameter point $P_{\ff}$ is generic in the sense of its coordinates;
that is,   the corresponding vector of coefficients belongs to a suitable dense Zariski open set of the parameter space.

Though the definition sounds a bit fluid, when working with such general set of forms, in practice, one deals with some property $\mathcal{P}$  that holds for them, meaning that there is a dense open set $U_{\mathcal{P}}$ (depending on $\mathcal{P}$) of the parameter space such that $\mathcal{P}$ holds for every set of forms whose corresponding point lies in $U_{\mathcal{P}}$.
 The stronger the given property, the smaller, and possibly the harder to describe, the corresponding dense open set in the parameter space, and the objective is often attained by indirect argument.

\begin{Proposition}\label{general_aci_satisfies_(a)}
	A general set of four forms of degree $d\geq 1$ in $R=\kk[x,y,z]$
	is $(d,4,e)$-extremal, for arbitrary $e\geq 1$ and generates an ideal of $R$ with socle degree $2d-2$.
\end{Proposition}
\demo
By Lemma~\ref{extremal_properties}, it suffices to show the existence of one set of four forms of degree $d$ satisfying the statement.

 Let then $f_1,f_2,f_3\in R=\kk[x,y,z]$ be a regular sequence of forms of degree $d$ and set $J:=(f_1,f_2,f_3)$. The minimal graded free resolution of $R/J$ is given by the Koszul complex
$$0\to R(-3d)\to R(-2d)^3\to R(-d)^3\to R\to R/J\to0,$$
from which one readily gets the Hilbert series of $R/J$:
$$H_{R/J}(t)=\frac{1-3t^d+3t^{2d}-t^{3d}}{(1-t)^3}.$$
Therefore,  the coefficients $a_i$ of  $H_{R/J}(t)$ are:

\begin{equation}\label{coeff_ci}
	a_i=\left\{\begin{array}{cccc}
		\dim_{\kk} R_i& \mbox{if $0\leq i\leq d-1$}\\
		\dim_{\kk} R_i-3\dim_{\kk} R_{i-d},&\mbox{if $d\leq i\leq 2d-1$}\\
		\dim_{\kk} R_{i}-3\dim_{\kk} R_{i-d}+3\dim_{\kk}R_{i-2d},&\mbox{if $2d\leq i\leq 3d-3$.}
	\end{array}\right.
\end{equation}
Now, if $\{f_1,f_2,f_3\}$ is a  general set then  $R/J$  has the Strong Lefschetz property (\cite{Stan}), that is, there is a  linear form $L\in R$  such that for $f:=L^d$  the multiplication map  $[R/J]_{i-d}\to [R/J]_{i}$ by $f$ has maximal rank for all $i\geq d$. 
Consequently, setting $I:=(f_1,f_2,f_3,f)$, since the image of this map is the vector space
$(fR_{i-d}, J_i)/J_i = I_i/J_i$, then 
$$\dim_{\kk} I_i/J_i=\min\{\dim_{\kk} [R/J]_{i-d},  \dim_{\kk} [R/J]_i\},$$
hence
\begin{eqnarray*}
\dim_{\kk} [R/I]_i &=& \dim_{\kk} [R/J]_i - \min\{\dim_{\kk} [R/J]_{i-d}, \dim_{\kk} [R/J]_i\}\\
&=& \max\{\dim_{\kk} [R/J]_i-\dim_{\kk} [R/J]_{i-d},\, 0\}
\end{eqnarray*}
for all $i\geq d$.

Thus, by \eqref{coeff_ci} and an obvious calculation, we have:
\begin{equation*}
	\dim_{\kk} [R/I]_{i}=\left\{\begin{array}{cccc}
		\dim_{\kk}R_i& \mbox{if $0\leq i\leq d-1$}\\
		\max\{\dim_{\kk}R_i-4\dim_{\kk}R_{i-d},\, 0\},&\mbox{if $d\leq i\leq 2d-1$}\\
		\max\{\dim_{\kk}R_i-4\dim_{\kk}R_{i-d}+3\dim_{\kk}R_{i-2d},\, 0\},&\mbox{if $2d\leq i\leq 3d-3$.}
	\end{array}\right.
\end{equation*}
Therefore, it follows that 
$$\dim_{\kk}[R/I]_i = \max\{\dim_{\kk} R_i - 4 \dim_{\kk} R_{i-d},\, 0\},$$
for arbitrary $i\geq d$, as was to be shown.
In particular, one has
$$\dim_{\kk} [R/I]_{i}=\dim_{\kk}R_{i}-4\dim_{\kk}R_{i-d}>0,\quad \mbox{for $d\leq i\leq 2d-2$}$$ 
and $\dim_{\kk} [R/I]_{2d-1}=0$, thus showing that the socle degree is $2d-2$.
\qed

\begin{Remark}\rm 
	The value of the socle degree, as well as the graded free resolution of $I$, has been previously established in \cite{Migl-Roig}.
\end{Remark}

As a last issue in this section we highlight a connection between Froberg's conjecture in arbitrary dimension and the concept of a $(d,r,e)$-extremal set of forms, in the way of relating it to \cite[Question 2.5]{Fro-Lund}.

Recall that, according to Fr\"oberg, one denotes by $((1-t^d)^r/(1-t)^n)_+$ the initial positive segment of the power series $(1-t^d)^r/(1-t)^n$.

\begin{Proposition}
	Let $R=\kk[x_1,\ldots,x_n]$ and $r\geq n$ and $d\geq 1$ be integers.
	Then
	\begin{equation*}
		\left(\frac{(1-t^d)^r}{(1-t)^n}\right)_+=\sum_{j=0}^{j_0-1}\left(\sum_{i=0}^{r}(-1)^{i}{r\choose i}\dim_{\kk}R_{j-id}\right)t^j
	\end{equation*}
	where $j_0$ is the least of the integers $d\leq j\leq dr-n$ satisfying  $\sum_{i=0}^{r}(-1)^{i}{r\choose i}\dim_{\kk}R_{j-id}\leq 0.$ 
\end{Proposition}	
\demo
The arbitrary term of the series $(1-t^d)^r/(1-t)^n$, is certainly well-known.
We recall how to get its expression.
Clearly, $1/(1-t)^n=\sum_{u=0}^{\infty}\dim_{\kk} R_{u}t^u.$ 
Thus,
\begin{eqnarray*}
\frac{(1-t^d)^r}{(1-t)^n}&=& \left(\sum_{i=0}^{r}(-1)^{i}{r\choose i}t^{di}\right)\left(\sum_{u=0}^{\infty}\dim_{\kk} R_{u}t^u\right)  \nonumber\\
	&=& \sum_{i=0}^{r}\sum_{u=0}^{\infty}(-1)^{i}{r\choose i}\dim_{\kk}R_{u}t^{u+id}\nonumber\\
	&=& \sum_{j=0}^{\infty} \left(\sum_{\substack{u\geq 0,\,0\leq i\leq r\\u+id=j}}(-1)^{i}{r\choose i}\dim_{\kk}R_{u}\right)t^{j} \nonumber\\
	&=&\sum_{j=0}^{\infty} \left(\sum_{i=0}^{r}(-1)^{i}{r\choose i}\dim_{\kk}R_{j-id}\right)t^{j}.
\end{eqnarray*}
From this and the definition of $((1-t^d)^r/(1-t)^n)_+$ the stated expression follows suit.
\qed

\begin{Remark}\rm 
The above clearly translates \cite[Question 2.5]{Fro-Lund} into the quest for an explicit expression of the integer $j_0$ in terms of $n, d, r$.
Also, we see that an affirmative answer to  Fr\"oberg's conjecture implies that  a general set of  $r\geq n$ forms $\ff=\{f_1,\ldots,f_r\}$  of degree $d$ is  $(d,r,e)$-extremal for every $0\leq e\leq \min\{d-1,j_0-d\}.$ 
\end{Remark}

\subsection{Gorenstein ideals and general sets of  forms}

We now draw some consequences of the previous section in dimension three.
First, one has the following constraint for Gorenstein ideals in this dimension.

\begin{Proposition}\label{linear_Gorenstein_is_not_general}
	Let $r\geq 5$ be an odd integer and $d=(r-1)/2.$ 
If there is  a Gorenstein ideal $I\subset R=\kk[x,y,z]$ generated by a general set of $r$ forms of degree $d$, then  $d=2$, i. e., $r=5$.
\end{Proposition}
\demo Let $I$ denote such an ideal.
By the symmetry of the Hilbert function of $R/I$ we have:
$$h_{R/I}(d+1)=h_{R/I}(d-3)={d-1\choose 2}.$$ 
Then, since $I$ is generated by a general set of forms, drawing upon Lemma~\ref{extremal_properties} (iii), we get
$${d-1\choose 2}=h_{R/I}(d+1)=
\max\left\{{d+3 \choose 2}-3r,\, 0\right\}
=\max\left\{{d+3\choose 2}-3(2d+1),\, 0\right\}.$$
A direct calculation shows that the only possibility for $d$ a positive integer is that  ${d-1\choose 2}=0$, hence $d=2$.
\qed

The main theorem is now a simple consequence: 

\begin{Theorem}\label{Old_Theorem 2.7} Let  $I\subset R=\kk[x,y,z]$ be a Gorenstein ideal generated by a general set of $r\geq 5$  forms of degree $d\geq 2$. Then $r=5$ and $d=2$.
	\end{Theorem}
\demo Let $J\subset I$ be generated by $n+1=4$ of the forms.
	By Proposition~\ref{general_aci_satisfies_(a)}, the socle degree of $R/J$ is $2d-2$. Then the socle degree of $R/I$ is at most $2d-2$.
	But, since $I$ is Gorenstein, the socle degree of $R/I$ is $2d+d'-3$ (Proposition~\ref{socle_degree_general}). Therefore, $d'=1$.
	Now apply Proposition~\ref{linear_Gorenstein_is_not_general}.
	\qed

\begin{Remark}\rm
An alternative proof of the above result is available by cooking up a couple of explicit calculations. Since these are rather lengthy and distracting, we have decided to defer them to the Appendix (Theorem~\ref{second_proof}).
\end{Remark}

More generally, we state the following conjecture.

\begin{Conjecture}
	Let $I\subset \kk[x_1,\ldots,x_n]$ be an ideal generated by a general set of $r\geq n+2$  forms of degree $d\geq 2$. Then $I$ is Gorenstein if and only if $d=2$ and $r= {{n+1}\choose 2}-1$.
\end{Conjecture}

We next proceed  towards the converse implication of this conjecture, by focusing on quadric generators.
To carry on,  we elaborate on some further notation.
Given integers $n\geq 3, d\geq 2$ and $r\geq 5$, consider the following ${{d+n-1}\choose d}\times r$ generic matrix
\begin{equation}\label{generic_coef}
\YY:=\YY_{n,d,r}=\left(\begin{array}{ccccccc}Y_{d,\ldots,0}^{(1)}&\cdots& Y_{d,\ldots,0}^{(r)}\\
\vdots&\ddots&\vdots\\
Y_{0,\ldots,d}^{(1)}&\cdots& Y_{0,\ldots,d}^{(r)}\end{array}\right),
\end{equation}
whose entries along the columns can be thought of as the respective coordinates of $r$ copies of $\P^{{{d+n-1}\choose d}-1}$.

One way of thinking about $r$  forms $\ff=\{f_1,\ldots,f_r\}\subset R=\kk[\xx]=\kk[x_1,\ldots,x_n]$ of degree $d\geq 1$ is as the matrix product $[\xx_d]\cdot \YY(P_{\ff})$, where $\xx_d$ is the list of monomials of degree $d$ in $\xx$ and $\YY(P_{\ff})$ denotes the  ${{d+n-1}\choose d}\times r$ matrix $\YY$ evaluated at $P_{\ff}$ (hence, with entries in $\kk$).

Throughout, set $\fm=(\xx)$.
Our main technical result is as follows.

\begin{Proposition}\label{Soc1bis}
	Let $I\subset R=\kk[x_1,\ldots,x_n]$ be an ideal generated by a general set of $n\leq r\leq N-1$ quadrics, where $N={{n+1} \choose 2}$. Then the socle of $R/I$ has no nonzero linear forms.
\end{Proposition}

\demo A quadric in $\kk[x_1,\ldots, x_n]$ depends on $N={{n+1}\choose 2}$ coefficients.
We will accordingly denote the corresponding indices by $\{1,1\},\{1,2\}, \ldots, \{n-1,n\}, \{n,n\}$.
Thus, the matrix $\YY_{n,2,r}$ takes the shape
\begin{equation*}
\left(\begin{array}{ccccc}
Y^{(1)}_{1,1}&\cdots&Y^{(r)}_{1,1}\\
\vdots&\ddots&\vdots\\
Y^{(1)}_{n,n}&\cdots&Y^{(r)}_{n,n}
\end{array}\right),
\end{equation*}
where delimiters have been omitted.
Let $\ff=\{f_1,\ldots,f_r\}\subset R$ be a general set of quadrics generating $I$. For simplicity,  set $P:=P_{\ff}.$ 
As explained above, one has the following matrix equality
$$[\ff]=[\xx_2]\cdot\YY_{n,2,r}(P).$$
The $N\times r$ matrix $\YY_{n,2,r}$ can be decomposed  in two vertical blocks as
$$\YY_{n,2,r}=\left(\begin{array}{c}B\\L\end{array}\right)$$
where
$$B=\left(\begin{array}{ccccc}
Y^{(1)}_{1,1}&\cdots&Y^{(r)}_{1,1}\\
\vdots&\ddots&\vdots\\
Y^{(1)}_{u,v}&\cdots&Y^{(r)}_{u,v}
\end{array}\right)_{r\times r}\quad \mbox{and}\quad L=\left(\begin{array}{ccccc}
Y^{(1)}_{u',v'}&\cdots&Y^{(r)}_{u',v'}\\
\vdots&\ddots&\vdots\\
Y^{(1)}_{n,n}&\cdots&Y^{(r)}_{n,n}
\end{array}\right)_{(N-r)\times r},$$
for suitables $1\leq u\leq v\leq n$ and $1\leq u'\leq v'\leq n$. 
Thus, 
\begin{equation}\label{comb_ff}
[\ff]\cdot{\rm cof}(B)(P)= [\xx_2]\cdot \left(\begin{array}{c}\Delta(P)\mathbb{I}_r\\L(P)\cdot{\rm cof}(B)(P)\end{array}\right)
\end{equation}
where ${\rm cof}(B)$  is the matrix of cofactors of $B$, $\mathbb{I}_r$ is the identity matrix of order $r$ and $\Delta=\det B$.
 Write, say,
\begin{equation}\label{LadjB}
L\cdot{\rm cof}(B)=
\left(\begin{array}{cccccccc}
\mathfrak{g}_{1,1}^{u',v'}&\cdots&\mathfrak{g}_{u,v}^{u',v'}\\
\vdots&\ddots&\vdots\\
\mathfrak{g}_{1,1}^{n,n}&\cdots&\mathfrak{g}_{u,v}^{n,n}
\end{array}\right).
\end{equation}

Since  $\ff$ is a general set, then $\Delta(P)\neq 0.$ Thus, by \eqref{comb_ff}:

$$I=\left(\Delta(P) x_ix_j+\sum_{(u',v')\leq (s,t)\leq (n,n)}\mathfrak{g}_{i,j}^{s,t}(P)x_sx_t \,|\, (1,1)\leq (i,j)\leq (u,v) \right).$$
Let $\ell=a_1x_1+\cdots+a_nx_n$ be a linear form in the socle of $R/I$. For any given $1\leq k\leq n$, one has
{\small\begin{equation}\label{mxl}
	a_1x_1x_k+\cdots+a_nx_nx_k=\sum_{(1,1)\leq (i,j)\leq (u,v)}\alpha_{i,j}\left(\Delta(P) x_ix_j+\sum_{(u',v')\leq (s,t)\leq (n,n)}\mathfrak{g}_{i,j}^{s,t}(P)x_sx_t\right).
	\end{equation}}
Comparing coefficients we get a  linear system $$\mathcal{A}_r(P)\cdot \aa=0$$ 
where $\aa$ is the transpose of the matrix $(a_1\cdots a_n)$ and  $\mathcal{A}_r(P)$ is an $(N-r)n\times n$ matrix whose  entries belong to the set 
$$\{0\}\cup\{\mathfrak{g}_{i,j}^{s,t}(P)\,:\; (1,1)\leq (i,j)\leq (u,v),\, (u',v')\leq (s,t)\leq (n,n)  \}\cup \{\Delta(P)\}.$$

We claim that $\rk(\mathcal{A}(P)_ r)=n$, which says that the assumed linear form is zero. We divide the proof in two cases:

\medskip

{\sc Case 1:} $r=N-1.$

In this case, the relation\eqref{mxl} has the following format for $1\leq t\leq n$:
\begin{equation}\label{linearinsocletimes}
a_1x_1x_t+\cdots+a_nx_nx_t=\sum_{\substack{1\leq i\leq j\leq n\\(i,j)\neq (n,n)}}\alpha_{i,j}(\Delta(P) x_ix_j-\mathfrak{g}_{i,j}^{n,n}(P)x_n^2).
\end{equation}
Comparing coefficients yields
\begin{equation}\label{linearsystem}
\mathcal{A}_{N-1}(P)=\left(\begin{array}{cccccc}\mathfrak{g}_{1,1}(P)&\cdots&\mathfrak{g}_{1,n-1}(P)&\mathfrak{g}_{1,n}(P)\\
\vdots&\ddots&\vdots&\vdots\\
\mathfrak{g}_{1,n-1}(P)&\cdots&\mathfrak{g}_{n-1,n-1}(P)&\mathfrak{g}_{n-1,n}(P)\\
\mathfrak{g}_{1,n}(P)&\cdots&\mathfrak{g}_{n-1,n}(P)&\Delta(P)\end{array}\right).
\end{equation}

Since the entries of  $\YY_{n,2,N-1}$ are mutually independent indeterminates, there is a $\kk$-homomorphism $\phi:\kk[\YY_{n,2,N-1}]\to \kk[L]$ mapping $B$ to the identity matrix and fixing the  entries of $L.$ 
Thus, by \eqref{LadjB}
$$L = [\phi(\mathfrak{g}_{1,1})\cdots\phi(\mathfrak{g}_{1,n})\,\,\phi(\mathfrak{g}_{2,2})\cdots\phi(\mathfrak{g}_{2,n})\cdots \phi(\mathfrak{g}_{n-1,n-1})\,\phi(\mathfrak{g}_{n-1,n})].$$
In particular, the determinant
$$\mathbb{D}:=\det\left(\begin{array}{cccccc}\mathfrak{g}_{1,1}&\cdots&\mathfrak{g}_{1,n-1}&\mathfrak{g}_{1,n}\\
\vdots&\ddots&\vdots&\vdots\\
\mathfrak{g}_{1,n-1}&\cdots&\mathfrak{g}_{n-1,n-1}&\mathfrak{g}_{n-1,n}\\
\mathfrak{g}_{1,n}&\cdots&\mathfrak{g}_{n-1,n}&\Delta\end{array}\right)$$
specializes to
$$\phi(\mathbb{D})=\det\left(\begin{array}{cccccc}Y_{n,n}^{(1)}&\cdots&Y_{n,n}^{(n-1)}&Y_{n,n}^{(n)}\\
\vdots&\ddots&\vdots&\vdots\\
Y_{n,n}^{(n-1)}&\cdots&Y_{n,n}^{(N-1)}&Y_{n,n}^{(N)}\\
Y_{n,n}^{(n)}&\cdots&Y_{n,n}^{(N)}&1\end{array}\right).$$
The latter does not vanish as it is the sum of two forms in two different degrees, none of which vanishes. 
Hence,  $\mathbb{D}\neq 0$ as well. 

\medskip

{\bf Case 2:} $r<N-1.$

It is enough to find a point $P'$ such that $\rk(\mathcal{A}_r(P'))=n.$ For this, consider a general  set of quadrics $\ff'=\{f'_1,\ldots,f'_{N-1}\}$. By the previous case, as applied to $J:=(f'_1,\ldots,f'_{N-1}),$ the socle  of $R/J$ has no linear forms.
Since $r\geq n$ and $\{f'_1,\ldots,f'_{N-1}\}$ is general, we can assume that the $r$ first elements $f'_1,\ldots,f'_r$ of $\ff'$ are such that $J'=(f'_1,\ldots,f'_r)$ is $\fm$-primary and $\Delta(P')\neq 0$, where $P'\in (\mathbb{P}^{N-1})^r$ is the point corresponding to the set $\{f'_1,\ldots,f'_r\}.$ 

If $\rk (\mathcal{A}_r(P'))<n$ then $\mathcal{A}_r(P')\cdot\aa=0$ has a nonzero solution $\aa$, so ${\rm Soc}(R/J')_1\neq 0.$ But, this is nonsense since  ${\rm Soc}(R/J')_1\subset {\rm Soc}(R/J)_1$ because  $J'\subset J$ live in degree $2$.
\qed

\medskip


The next result follows suit.

\begin{Theorem}\label{socle_always_in_degree2}
	Let $I\subset R=\kk[x_1,\ldots,x_n]$ be an ideal generated by a general set of $\frac{(n+2)(n+1)}{6}\leq r\leq N-1$  quadrics, where $N={{n+1}\choose 2}$. Then ${\rm Soc}(R/I)=\kk(-2)^{N-r}.$
	In particular, $R/I$ is Gorenstein if and only if $r=N-1$.
\end{Theorem}

\demo We have 
\begin{eqnarray}
\dim_{\kk} [R/I]_3&=&\max\{\dim_{\kk} R_3-rn,\, 0\}\quad \mbox{(by Lemma~\ref{extremal_properties} (iii)) }
	\nonumber\\
&=& 0 \quad \quad\quad\quad\quad\quad\quad\quad\quad\,\,\,(\mbox{because $(n+2)(n+1)/6\leq r\leq N-1$}).\nonumber
\end{eqnarray} 
 Thus, $R/I=\kk\oplus \kk(-1)^{n}\oplus \kk(-2)^{N-r}$
 as a $\kk$-vector space.
Hence, by Proposition~\ref{Soc1bis}, 
$${\rm Soc}(R/I)= \kk(-2)^{N-r}.$$
In particular, $\dim_{\kk} {\rm Soc}(R/I)=1$ if and only if $r=N-1,$ that is, $R/I$ is an Artinian Gorenstein algebra if and only if  $r=N-1.$
\qed




\begin{Remark}\rm
A simpler proof is available in the case where $n=3$.
Namely, by the argument in the proof of \cite[Proposition 2.3]{abc},  if the socle of $I$ contains a linear form then $\mu(I^2)<15$. However, if $I$  is generated by a general set of five  forms, we have $\mu(I^2)=15$, thus giving a contradiction. 
This argument breaks down for $n\geq 4$, whereas the theorem fixes it for $n\leq 5$.
\end{Remark}

\section{On the $(x_1^{m},\ldots,x_n^m)$-colon problem}

It is known (see \cite[Proposition 1.3]{BE0}) that any homogeneous Gorenstein ideal of codimension $n$ in $\kk[\xx]=\kk[x_1,\ldots,x_n]$ can be obtained as a colon ideal $(x_1^m,\ldots,x_n^m):\mathfrak{f}$, for some integer $m \geq 1$ and some form $\mathfrak{f}$. 
In this section we deal with some of the main questions regarding this representation.

It is first established under which condition the form $\mathfrak{f}$ is uniquely determined and what is its degree  in terms of $m$ and the socle degree of $I$.
Then we prove that $\mathfrak{f}$ can be retrieved from $I$ by taking the so-called (socle-like) Newton dual of a minimal generator of the Macaulay inverse of $I$.

Then we give conditions under which the Gorenstein ideal $I$ is equigenerated in terms of the exponent $m$ and the form $\mathfrak{f}$.
We solve this problem in the case where $I$ has linear resolution
$$0 \rar R(-2d -n + 2) \rar R(-d - n + 2) ^{b_{n-1}} \rar\cdots\rar R(-d - 1)^{b_2} \rar R(-d)^{b_1} \rar  R$$
where $b_1=\mu(I)$.

These questions will be subsumed under the designation {\em the colon problem}, to avoid `link' which has already many uses.
For convenience, call $\mathfrak{f}$ a {\em directrix form} (of $I$) associated to the regular sequence $\{x_1^{m},\ldots,x_n^m\}$.

\subsection{Macaulay inverse system versus Newton duality}\label{Mac_inverses}


We briefly recall  some main features of the Macaulay inverse system.
For recent accounts of this classical theme see, e. g.,  \cite{Elias} ,  \cite{EliasRossi} and \cite{Kustin}.

Let $V$ be a vector space of dimension $n$ over a field $\kk$ and let $x_1,\ldots,x_n$ be a basis for $V.$  Let $R={\rm Sym}_{\kk}(V)=\kk[x_1,\ldots,x_n]$ be the standard graded polynomial ring in $n$ variables over $\kk.$ Set $y_1,\ldots,y_n$ for the dual basis on  $V^{\ast}=\Hom_{\kk}(V,\kk)$ and consider the divided power ring $$D_{\kk}(V^{\ast})=\displaystyle\bigoplus_{i\geq 0}\Hom(R_i,\kk)=\kk_{\rm DV}[y_1,\ldots,y_n].$$ 
In particular, $\{\yy^{[\alpha]}\,|\,\alpha\in\mathbb{N}^n\,\mbox{and}\,|\alpha|=j\}$ is the dual basis of $\{\xx^{\alpha}\,|\,\alpha\in\mathbb{N}^n\,\mbox{and}\,|\alpha|=j\}$ on $D_{\kk}(V^{\ast})_j=\Hom(R_j,\kk).$ If $\alpha\in\mathbb{Z}^n$ then we set $\yy^{[\alpha]}=0$ if some component of $\alpha$ is negative.
Make $ D_{\kk}(V^{\ast})$ into a module over $R$ through the following action 
$$R\times D_{\kk}(V^{\ast})\to D_{\kk}(V^{\ast}),\quad (f=\sum_{\alpha}a_{\alpha}\xx^{\alpha},F=\sum_{\beta}b_{\beta}\yy^{\beta})\mapsto fF=\sum_{\alpha,\beta}a_{\alpha}b_{\beta}\yy^{[\beta-\alpha]}.$$

For a homogeneous ideal $I\subset R$ and an $R$-submodule $M\subset D_{\kk}(V^{\ast})$ one defines:

$${\rm Ann}(I):=\{g\in D_{\kk}(V^{\ast})\,|\,Ig=0 \}\quad \mbox{and}\quad {\rm Ann}(M):=\{f\in R\,|\, fM=0\}.$$

Then ${\rm Ann}(I)$ is an $R$-submodule of $D_{\kk}(V^{\ast})$, while ${\rm Ann}(M)$ is an ideal of $R.$ The $R$-module ${\rm Ann}(I)$ is called the {\it Macaulay inverse} (system) of $I.$

The main basic result regarding this construction is due to Macaulay (\cite{Macaulay}).
In the present  language it can be stated in the following form:

\smallskip

\begin{Theorem}\label{Macaulay_main}{\rm  (Macaulay Duality, (\cite[Theorem 1.4]{Kustin})}
	There exists a one-to-one correspondence between the set of nonzero homogeneous height $n$ Gorenstein ideals of $R$ and the set of nonzero homogeneous cyclic submodules of $D_{\kk}(V^{\ast})$ given by $I\mapsto {\rm Ann}(I)$ with inverse $M\mapsto {\rm Ann}(M).$ Moreover, the socle degree of $R/I$ is equal to the degree of a homogeneous generator of ${\rm Ann}(I).$
\end{Theorem} 

The Macaulay--Matlis duality meets yet another version in terms of the Newton polyhedron nature of the homogeneous forms involved so far.

For this, recall the notion of the Newton (complementary) dual of a form $f\in\kk[\xx]=\kk[x_1,\ldots,x_n]$ in a polynomial ring over a field $\kk$, as introduced in  \cite{CostaSi}, and \cite{DoriaSi}. 
Namely, start out with  the log matrix $A$ of the constituent monomials of $f$ (i.e, the nonzero terms of $f$).
This is the matrix whose columns are the exponents vectors of the nonzero terms of $f$ in, say, the lexicographic ordering. 
It is denoted $\mathcal{N}(f)$.
Then, the  Newton dual log matrix (or simply the  {\em Newton dual matrix}) of the Newton log matrix $\mathcal{N}(f)=(a_{i,j})$ is the matrix $\widehat{\mathcal{N}(f)}=(\alpha_i - a_{i,j}),$
where $\alpha_i = \max_j \{a_{i,j}\}$, with $1\leq i\leq n$
and $j$ indexes the set of all nonzero terms of $f$.

In other words, denoting
${\boldsymbol\alpha}:=(\alpha_1\cdots\alpha_n)^t$, one has
$$
\widehat{\mathcal{N}(f)}=\left[\,{\boldsymbol\alpha}\,| \cdots
|\,{\boldsymbol\alpha}\,\right]_{(n+1) \times r}
- \mathcal{N}(f),
$$
where $r$ denotes the number of nonzero terms of $f$.
The vector ${\boldsymbol\alpha}$ is called the {\em directrix vector} of $\mathcal{N}(f)$ (or of $f$ by abuse).

We note that taking the Newton dual is a true duality upon forms not admitting monomial factor, in the sense that, for such a form $f$, $\widehat{\widehat{\mathcal{N}(f)}}=\mathcal{N}(f)$ holds.

We define the {\em Newton dual} of $f$ to be the form $\hat{f}$ whose terms are the ordered monomials obtained form $\widehat{\mathcal{N}(f)}$ affected by the same coefficients as in $f$.



Our next result asserts that directrix forms and Macaulay inverse generators obey a duality in terms of the above Newton dual.
Given a  directrix form $\mathfrak{f}$ associated to the regular sequence $\{x_1^{m},\ldots,x_n^{m}\}$ -- i.e., $(x_1^{m},\ldots,x_n^{m}):\mathfrak{f}=I$ -- it will typically admit monomial terms belonging to the ideal $(x_1^{m},\ldots,x_n^{m})$.
In order to fix this inconvenient, we redefine the {\em socle-like Newton dual} of such directrix form by taking as directrix vector $\nu:=(m-1\, \cdots\, m-1)^t$.

\begin{Proposition}\label{NewtondualxMinverse}
	Let $I\subset R=\kk[\xx]=\kk[x_1,\ldots,x_n]$ be a homogeneous codimension $n$ Gorenstein ideal with socle degree $s$.
Given an integer $m\geq 1$, suppose that $I$ admits a directrix form $\mathfrak{f}$ associated to the regular sequence $\{x_1^{m},\ldots,x_n^{m}\}$.
	Then:
	\begin{enumerate}
		\item[{\rm (i)}] $\mathfrak{f}$ is a degree $n(m-1)-s$ form uniquely determined, up to a scalar coefficient, by the condition that no nonzero term of $\mathfrak{f}$ belongs to the ideal $(x_1^{m},\ldots,x_n^{m})$.
		\item[{\rm (ii)}] The socle-like Newton dual of $\mathfrak{f}$ is a minimal generator of the Macaulay inverse to $I$ {\rm (}having dual degree $s${\rm )}, and its socle-like Newton dual retrieves $\mathfrak{f}$.  
	\end{enumerate}  
\end{Proposition}
\demo  Suppose $\mathfrak{f}=\displaystyle\sum_{|\alpha|=\deg \mathfrak{f}} a_{\alpha}\xx^{\alpha}.$ Then, the socle-like Newton dual of $\mathfrak{f}$ is $\hat{\mathfrak{f}}=\displaystyle\sum_{\alpha} a_{\alpha}\yy^{\hat{\alpha}},$ where $\hat{\alpha}:=\nu-\alpha$ (in particular, $\deg \hat{\mathfrak{f}}=n(m-1)-\deg\mathfrak{f}$).  Given a homogeneous polynomial $h=\displaystyle\sum_{|\beta|=\deg h}b_{\beta}\xx^{\beta}\in R$ one has:
$$h\mathfrak{f}=\sum_{|\gamma|=\deg\mathfrak{f}+\deg h}\left(\sum_{\alpha+\beta=\gamma} a_{\alpha}b_{\beta}\right)\xx^{\gamma}\quad\mbox{and}\quad h\hat{\mathfrak{f}}=\sum_{|\gamma|=\deg\mathfrak{f}+\deg h}\left(\sum_{\hat{\alpha}-\beta=\hat{\gamma}} a_{\alpha}b_{\beta}\right)\yy^{\hat{\gamma}}$$
with $\hat{\gamma}=\nu-\gamma.$ In particular, for every $\gamma,$ the coefficient of $\xx^{\gamma}$ as a term in $h\mathfrak{f}$ is equal to the coefficient of $\yy^{\hat{\gamma}}$ as a term of  $h\hat{\mathfrak{f}}$. Moreover, the $i$th coordinate of $\gamma$ is larger than $s$ if and only if the $i$th coordinate of $\hat{\gamma}$ is negative. Thus,
\begin{eqnarray}
h\in (x_1^{m},\ldots,x_n^{m}):\mathfrak{f} &\Leftrightarrow& \mbox{for every} \sum_{\alpha+\beta=\gamma} a_{\alpha}b_{\beta}\neq 0, \gamma\,\mbox{has a coordinate larger than}\,m\nonumber\\
&\Leftrightarrow& \mbox{for every} \sum_{\hat{\alpha}-\beta=\hat{\gamma}} a_{\alpha}b_{\beta}\neq 0, \hat{\gamma}\,\mbox{has a negative coordinate}\\
&\Leftrightarrow& h\hat{\mathfrak{f}}=0 \Leftrightarrow h\in{\rm Ann}(\hat{\mathfrak{f}}).\nonumber
\end{eqnarray}
Therefore, $I={\rm Ann}(\hat{\mathfrak{f}}),$ that is, $\hat{\mathfrak{f}}$ is a minimal generator of the Macaulay inverse to $I$. By construction, one has $\deg \hat{\mathfrak{f}}=n(m-1)-\deg {\mathfrak{f}}.$ On the other hand, it is well known that the degree of a minimal generator of the Macaulay inverse to $I$ is the socle degree of $I$, i.e., $\deg\hat{\mathfrak{f}} =s$. 
Therefore, $\deg\mathfrak{f} =n(m-1)-s.$ Since $\hat{\mathfrak{f}}$ is uniquely determined, up to a scalar coefficient, the form  $\mathfrak{f}$ is  uniquely determined as well, up to a scalar coefficient, by the condition that no nonzero term of $\mathfrak{f}$ belongs to the ideal $(x_1^{m},\ldots,x_n^{m})$. Thus, assertion (i) follows. 

Assertion (ii) follows from the above.
\qed

\begin{Remark}\label{from_vars_to_linear_forms}\rm
Item (i) of Proposition~\ref{NewtondualxMinverse} is stable under a change of coordinates.
In other words, it holds true replacing the sequence $\{x_1^m,\ldots, x_n^m\}$ by a sequence $\{\ell_1^m,\ldots, \ell_n^m\}$, where $\{\ell_1,\ldots,\ell_n\}$ are independent linear forms.
Thus, if $I$ is a homogeneous codimension $n$ Gorenstein ideal such that $(\ell_1^m,\ldots, \ell_n^m):\mathfrak{f}=I$, for some form $\mathfrak{f}\in R$, then $\mathfrak{f}$ is uniquely determined, up to a scalar coefficient, by the condition that no nonzero term  of $\mathfrak{f}$,
 written as a polynomial in $\ell_1,\ldots,\ell_n,$ belongs to the ideal $(\ell_1^{m},\ldots,\ell_n^{m})$.
\end{Remark}

\subsection{The case of linear resolution}

In this section we  characterize when  $I$ is an equigenerated codimension $n$ Gorenstein ideal with linear resolution in terms of the exponent $m$ and the form $\mathfrak{f}\in R=\kk[\xx]=\kk[x_1,\ldots,x_n].$ The preliminaries remain valid in arbitrary characteristic, but  characteristic zero is called upon in item (ii) of Proposition~\ref{sum_power_case} below.

Let $e,e'$ be positive integers and let $\mathfrak{f}=\sum_{|\alpha|=e}a_{\alpha}\xx^{\alpha}\in R_{e}$ and $g=\sum_{|\beta|=e'}b_{\beta}\xx^{\beta}\in R_{e'}$ be forms.
Given an integer $m\geq 1$, write
\begin{equation}\label{gf}
g\mathfrak{f}=\sum_{\xx^{\gamma}\notin (x_1^m,\ldots,x_n^{m})}\left(\sum_{\alpha+\beta=\gamma} a_{\alpha}b_{\beta}\right)\xx^{\gamma}+\sum_{\xx^{\gamma}\in (x_1^m,\ldots,x_n^{m})}\left(\sum_{\alpha+\beta=\gamma} a_{\alpha}b_{\beta}\right)\xx^{\gamma},
\end{equation}
where $\gamma\in \mathbb{N}^n$ is a running $n$-tuple.
To this writing associate a matrix $\mathcal{M}_{e,e',m}$ whose rows are indexed by the $n$-tuples $\gamma$ such that $|\gamma|=e+e'$ and whose columns are indexed by the $n$-tuples $\beta$ such that $|\beta|=e'.$ The entries of the matrix are specified as follows:

$$\text{the $(\gamma,\beta)$-entry of $\mathcal{M}_{e,e',m}$}=
\left\{\begin{array}{cc} 0,&  \text{if some coordinate of}\, \gamma-\beta \,\mbox{is}\,<0 \\ a_{\alpha},&\text{if each coordinate of}\, \alpha=\gamma-\beta \,\mbox{is}\,\geq 0.
\end{array}\right.$$

In addition, let $\chi$ denote the row matrix $[\xx^{\gamma}]$ with the monomial entries $\xx^{\gamma}\notin (x_1^m,\ldots,x_n^{m})$, and let $\bb$  stand for the column matrix whose entries are the coefficients $b_{\beta}$ of $g.$

Then equality\eqref{gf} can be rewritten in the shape

\begin{equation}
g\mathfrak{f}=\chi\cdot \mathcal{M}_{e,e',m}\cdot\bb+ \sum_{\xx^{\gamma}\in (x_1^m,\ldots,x_n^{m})}\left(\sum_{\alpha+\beta=\gamma} a_{\alpha}b_{\beta}\right)\xx^{\gamma}.
\end{equation}
It is important to observe that the matrix $\mathcal{M}_{e,e',m}$ depends only on the integers $e,e'$ and $m$, and not on the details of $g$.

From this, it follows immediately:

\begin{Lemma}\label{rankM}
	$g\in I=(x_1^m,\ldots,x_n^m):\mathfrak{f}$ if and only if $\mathcal{M}_{e,e',m}\cdot\bb=0.$ In particular, $I_{e'}=\{0\}$ if and only if $\rank \mathcal{M}_{e,e',m}={e'+n-1\choose n-1}.$
\end{Lemma}

To tie up the ends, consider the parameter map

$$R_e\to \mathbb{P}^{{e+n-1\choose n-1}-1},\quad \mathfrak{f}=\sum_{|\alpha|=e}a_{\alpha}\xx^{\alpha}\mapsto P_{\mathfrak{f}}=(a_{(e,\ldots,0)}:\cdots:a_{(0,\ldots,e)})$$
in the notation of Subsection~\ref{Parameters}.
Let $\{Y_{e,\ldots,0},\ldots,Y_{0,\ldots,e}\}$ denote the coordinates of $\mathbb{P}^{{e+n-1\choose n-1}-1}$ and let $\mathcal{MG}_{e,e',m}$ stand for the matrix whose entries  are obtained by replacing  each $a_{\alpha}$ in $\mathcal{M}_{e,e',m}$ by the corresponding $Y_{\alpha}.$

\begin{Theorem}\label{colon_problem_linear_case}
	Let $m\geq 1$ be an integer and let $\mathfrak{f}\in R=\kk[x_1,\ldots,x_n]$ be a form. The following are equivalent:
	\begin{enumerate}
		\item[{\rm(i)}] $I=(x_1^m,\ldots,x_n^m):\mathfrak{f}$ is a codimension $n$ equigenerated Gorenstein ideal with linear resolution.
		\item[{\rm(ii)}] The integer $s:=n(m-1)-\deg \mathfrak{f}$ is even  and $\rank \mathcal{M}_{\deg \mathfrak{f},s/2,m}={s/2+n-1\choose n-1}.$
		\item[{\rm(iii)}] The integer $s:=n(m-1)-\deg \mathfrak{f}$ is even and $P_{\mathfrak{f}}$ is a point in the Zariski open set $\mathbb{P}^{{e+n-1\choose n-1}-1}\setminus V(I_{k}(\mathcal{MG}_{e,s/2,m}))$, with $k:={s/2+n-1\choose n-1}$ and $e=\deg \mathfrak{f}$.
	\end{enumerate}
\end{Theorem}

\demo (i)$\Rightarrow$(ii) Suppose that $I$ is equigenerated in degree $d.$ Since $I$ has linear resolution then the socle degre of $I$ is $2d-2.$ Thus, by the Proposition~\ref{NewtondualxMinverse}, $s=2d-2.$ In particular, $s$ is an even integer. On the other hand, since $I$ is generated in degree $d$ then $I_{s/2}=I_{d-1}=\{0\}.$ Hence, by the Lemma \ref{rankM}, $\rank \mathcal{M}_{\deg\mathfrak{f},s/2,m}={s/2+n-1\choose n-1}.$

(ii)$\Rightarrow$(i)  We claim that $I$ is codimension $n$ Gorenstein ideal generated in degree $t=s/2+1.$ The ideal $I$ is Gorenstein of codimension $n$ because  it is the link of the homogeneous almost complete intersection $J=(x_1^m,\ldots,x_n^m,\mathfrak{f})$ with respect to the complete intersection of pure powers $(x_1^m,\ldots,x_n^m).$  By Proposition~\ref{NewtondualxMinverse}, the socle degree of $I$ is $s.$ Thus, $(R/I)_{2t-1}=(R/I)_{s+1}=\{0\}.$ On the other hand, since $\rank \mathcal{M}_{\deg\mathfrak{f},s/2,m}={s/2+n-1\choose n-1},$ then $I_{t-1}=I_{s/2}=\{0\}$ by Lemma~\ref{rankM}. 
Since $I$ is a codimension $n$ Gorenstein ideal and $(R/I)_{2t-1}=\{0\}$ and $I_{t-1}=0$ it follows from \cite[Proposition 1.8]{Kustin} that $I$ is generated in degree $t$ and has linear resolution.

(ii) $\Leftrightarrow$ (iii) This is a mere language transcription.
\qed

\begin{Remark}\rm The key point for proving the implication (i) $\Rightarrow$ (ii) is the use of \cite[Proposition 1.8]{Kustin}, which characterizes the  $\fm$-primary Gorenstein ideals with linear resolution through estimates for the initial degree and the socle degree. For other classes of equigenerated Gorenstein ideals the examples show that a similar characterization must take into account not only the initial degree and the socle degree.   For example,
 the ideals $(x^5,y^5,z^5):(x+y+z)^5$ and $(x^5,y^5,z^5):x^3y^2+y^3z^2+x^2z^3$ have the same initial degree and the same socle degree. However, the first ideal is equigenerated in degree 4 while the second is minimally generated in degree 4 and 5. Extending Theorem~\ref{colon_problem_linear_case} to other ideals should  include additional conditions.
	\end{Remark}

\smallskip

The question remains as to when the Zariski open set $\mathbb{P}^{{e+n-1\choose n-1}-1}\setminus V(I_{k}(\mathcal{MG}_{e,s/2,m}))$ is nonempty, where $k:={s/2+n-1\choose n-1}$ and $e=\deg \mathfrak{f}$. 
The next result determines all pair of integers $m,e\geq 1$, with even $s=n(m-1)-e$, for this to be the case when $\mathfrak{f}=(x_1+\cdots+x_n)^{e}$.

\begin{Proposition}\label{sum_power_case} {\rm (char$(\kk)=0$)}
	Let $m,e\geq 1$ integers such that $s=n(m-1)-e$ is even. Set $d:=s/2+1.$
	\item[{\rm(i)}] If $m< d$ then $\mathbb{P}^{{e+n-1\choose n-1}-1}\setminus V(I_{k}(\mathcal{MG}_{e,s/2,m}))=\emptyset.$
	\item[{\rm(ii)}] If $m\geq d$ then $I=(x_1^m,\ldots,x_n^m):(x_1+\cdots+x_n)^{e}$ is a codimension $n$ Gorenstein ideal  generated by forms of degree $d$ with linear resolution. In particular, $\mathbb{P}^{{e+n-1\choose n-1}-1}\setminus V(I_{k}(\mathcal{MG}_{e,s/2,m}))$ is a dense open set.
\end{Proposition}
\demo (i) We claim that there is no form $\mathfrak{f}$ of degree $e$ such that $I=(x_1^m,\ldots,x_n^m):\mathfrak{f}$ is a equigenerated codimension $n$ Gorenstein ideal with linear resolution. In fact, otherwise $I$ would be an ideal generated in degree $d$ with $(x_1^m,\ldots,x_n^m)\subset I$ --  an absurd. Hence, by Theorem~\ref{colon_problem_linear_case},  $\mathbb{P}^{{e+n-1\choose n-1}-1}\setminus V(I_{k}(\mathcal{MG}_{e,s/2,m}))=\emptyset.$

(ii)  We mimic the argument of \cite[Proposition 7.24]{Kustin}. Namely, by applying \cite[Proposition 1.8]{Kustin}, it is sufficient to show that $R_{2d-1}\subset I$ and $I_{d-1}=\{0\}.$ Clearly, $$R_{n(m-1)+1}\subset (x_1^m,\ldots,x_n^m).$$ Moreover, 
$$(x_1+\cdots+x_n)^{e}R_{2d-1}\subset R_{2d-1+e}=R_{n(m-1)+1}.$$ 
Hence, $R_{2d-1}\subset I.$ On the other hand, the initial degree of $I/(x_1^m,\ldots,x_n^m)$ is at least $d$ as a consequence of the Lefschetz like result of R. Stanley, as proved in \cite[Theorem 5]{RRR91}. Since $m\geq d$ by assumption then the initial degree of $I$ is at least $d$.
Therefore, $I_{d-1}=\{0\}$, as was to be shown. In particular, for $\mathfrak{f}=(x_1+\cdots+x_n)^{e}$ Theorem~\ref{colon_problem_linear_case} gives $P_{\mathfrak{f}}\in\mathbb{P}^{{e+n-1\choose n-1}-1}\setminus V(I_{k}(\mathcal{MG}_{e,s/2,m})).$
\qed

\begin{Remark}\rm
The only place where one needs characteristic zero above is in the use of \cite[Theorem 5]{RRR91} -- for a different proof of this typical charming result of characteristic zero, see \cite[Th\'eor\`eme, Appendix]{BuChSi}. It is reasonable to expect that the above proposition be valid in arbitrary characteristic.
\end{Remark}

\subsection{The pure power gap}

Let $ \boldsymbol\ell=\{\ell_1,\ldots,\ell_n\}\in R=\kk[x_1,\ldots,x_n]$ be a regular sequence of linear forms and let $I\subset R $  be a homogeneous codimension $n$ Gorenstein ideal with socle degree $s.$ Denote by $m(I,\boldsymbol\ell)$ the least index $m$ such that $\{\ell_1^m,\ldots,\ell_n^m\}\subset I.$ Since $R_{s+1}=I_{s+1},$  then $m(I,\boldsymbol\ell)\leq s+1.$ The {\it pure power gap} of $I$ with respect to the regular sequence $\boldsymbol\ell$  is $\mathfrak{g}(I,\boldsymbol\ell):=s+1-m(I,\boldsymbol\ell).$ The  {\it absolute pure power gap} of $I$ (or simply, the pure power gap of $I$) is $\mathfrak{g}(I):=s+1-\displaystyle\min_{\boldsymbol\ell}\{m(I,\boldsymbol\ell)\}.$

To start we have the following basic ring-theoretic result:

\begin{Lemma}\label{cond_iterado}
	Let $m_1,\ldots,m_n\geq 1$ be  integers and $\mathfrak{f}$ a form in $R=\kk[x_1,\ldots,x_n].$ Then $$(\ell_1^{m_1},\ldots,\ell_n^{m_n}):\mathfrak{f}=(\ell_1^{m_1},\ldots,\ell_i^{m_i+1},\ldots,\ell_n^{m_n}):\ell_i\mathfrak{f}$$
	for every $1\leq i\leq n.$ In particular, $$(\ell_1^{m_1},\ldots,\ell_n^{m_n}):\mathfrak{f}=(\ell_1^{m_1+k},\ldots,\ell_i^{m_i+k},\ldots,\ell_n^{m_n+k}):(\ell_1\cdots\ell_n)^k\mathfrak{f}$$
	for each $k\geq 0.$
\end{Lemma}
\demo One can assume that $i=1.$  The inclusion 
$(\ell_1^{m_1},\ldots,\ell_n^{m_n}):\mathfrak{f}\subset(\ell_1^{m_1+1},\ldots,\ell_n^{m_n}):\ell_1\mathfrak{f}$ is immediate. 
Thus, consider $h\in (\ell_1^{m_1+1},\ldots,\ell_n^{m_n}):\ell_1\mathfrak{f}.$ Then, 

$$\ell_1\mathfrak{f}h=p_1\ell_1^{m_1+1}+\cdots +p_n\ell_n^{m_n}$$ 
for certain $p_1,\ldots,p_n\in R.$  In particular, $\ell_1$ divide $p_2\ell_2^{m_2}+\cdots+p_n\ell_n^{m_n}.$ We can write 
$$p_i=\ell_1q_i+r_i,\quad \mbox{for each}\,2\leq i\leq n,$$
where $r_2,\ldots,r_n$ are polynomials in $\kk[\ell_2,\ldots,\ell_n].$ Thus, $$p_2\ell_2^{m_2}+\cdots+p_n\ell_n^{m_n}=q_2\ell_1\ell_2^{m_2}+\cdots+q_n\ell_1\ell_n^{m_n}+r_2\ell_2^{m_2}+\cdots+r_n\ell_n^{m_n}.$$
Since $\ell_1$ divides $p_2\ell_2^{m_2}+\cdots+p_n\ell_n^{m_n}$ and $r_2\ell_2^{m_2}+\cdots+r_n\ell_n^{m_n}\in \kk[\ell_2,\ldots,\ell_n]$ then 
$$p_2\ell_2^{m_2}+\cdots+p_n\ell_n^{m_n}=q_2\ell_1\ell_2^{m_2}+\cdots+q_n\ell_1\ell_n^{m_n}.$$
Thus, 
$$\mathfrak{f}h=p_1\ell_1^{m_1}+q_2\ell_2^{m_2}+\cdots +q_n\ell_n^{m_n},$$
that is, $h\in (\ell_1^{m_1},\ldots,\ell_n^{m_n}):\mathfrak{f}.$ 
Therefore, $(\ell_1^{m_1},\ldots,\ell_n^{m_n}):\mathfrak{f}=(\ell_1^{m_1+1},\ldots,\ell_n^{m_n}):\ell_1\mathfrak{f} $ as stated.
\qed

\smallskip

To see an application, recall from Remark~\ref{from_vars_to_linear_forms} that if  $I$ is a homogeneous codimension $n$ Gorenstein ideal such that $(\ell_1^m,\ldots, \ell_n^m):\mathfrak{f}=I$, where $\ell_1,\ldots, \ell_n$ are linear forms,  then $\mathfrak{f}$ is uniquely determined, up to a scalar coefficient, by the condition that no nonzero term  of $\mathfrak{f}$,
written as a polynomial in $\ell_1,\ldots,\ell_n,$ belongs to the ideal $(\ell_1^{m},\ldots,\ell_n^{m})$.

\begin{Proposition}\label{power_gap_characterization}  Let $I\subset R $  be a homogeneous codimension $n$ Gorenstein ideal with socle degree $s.$
Suppose that as above,  $\ell_1,\ldots, \ell_n$ are linear forms such that
  $I=(\ell_1^{s+1},\ldots,\ell_{n}^{s+1}):\mathfrak{f}$ with $\mathfrak{f}$ uniquely determined, up to a scalar coefficient, by the condition that no nonzero term of $\mathfrak{f}$ belongs to the ideal $(\ell_1^{s+1},\ldots,\ell_{n}^{s+1}).$ Then, $\mathfrak{g}(I,\boldsymbol\ell)$ is the largest index such that $(\ell_1\cdots\ell_n)^{\mathfrak{g}(I,\boldsymbol\ell)}$ divides $\mathfrak{f}.$
\end{Proposition}

\demo Denote $m_0:=m(I,\boldsymbol\ell)$ and $\mathfrak{g}:=\mathfrak{g}(I,\boldsymbol\ell).$ 
Then $(\ell_1^{m_0},\ldots,\ell_{n}^{m_0}):I$ is an almost complete intersection $J=(\ell_1^{m_0},\ldots,\ell_{n}^{m_0}, \mathfrak{f}_0)$, for some 
form $\mathfrak{f}_0\in R$. Since $R=\kk[\ell_1,\ldots,\ell_n]$, we can write 
$\mathfrak{f}_0$ as a polynomial in these linear forms and get rid of the terms belonging to the ideal $(\ell_1^{m_0},\ldots,\ell_{n}^{m_0})$.
This way, the latter is part of a minimal set of generators of $J$.
Therefore, $(\ell_1^{m_0},\ldots,\ell_{n}^{m_0}):J=(\ell_1^{m_0},\ldots,\ell_{n}^{m_0}):\mathfrak{f}_0$ is Gorenstein and $I=(\ell_1^{m_0},\ldots,\ell_{n}^{m_0}):\mathfrak{f}_0$.

 By  Lemma~\ref{cond_iterado} one has
\begin{eqnarray}
I=(\ell_1^{m_0},\ldots,\ell_{n}^{m_0}):\mathfrak{f}_0&=&(\ell_1^{m_0+\mathfrak{g}},\ldots,\ell_{n}^{m_0+\mathfrak{g}}):(\ell_1\cdots\ell_n)^{\mathfrak{g}}\mathfrak{f}_0\nonumber\\
&=&(\ell_1^{s+1},\ldots,\ell_{n}^{s+1}):(\ell_1\cdots\ell_n)^{\mathfrak{g}}\mathfrak{f}_0.
\end{eqnarray}
Consider  $$\mathfrak{f}_0=\sum_{|\alpha|=\deg\mathfrak{f}_0}a_{\alpha}\ell_1^{\alpha_1}\cdots\ell_n^{\alpha_n}.$$
Then, $$(\ell_1\cdots\ell_n)^{\mathfrak{g}}\mathfrak{f}_0=\sum_{|\alpha|=\deg\mathfrak{f}_0}a_{\alpha}\ell_1^{\alpha_1+\mathfrak{g}}\cdots\ell_n^{\alpha_n+\mathfrak{g}}.$$

For each nonzero $a_{\alpha},$ one has $\alpha_i\leq m_0-1$ for each $1\leq i\leq n.$ Hence, $\alpha_i+\mathfrak{g}\leq m_0+\mathfrak{g}-1=s$ for each $1\leq i\leq n.$ Thus, no nonzero term of $(\ell_1\cdots\ell_n)^{\mathfrak{g}}\mathfrak{f}_0$ belongs to the ideal $(\ell_1^{s+1},\ldots,\ell_{n}^{s+1}).$ Then, since $\mathfrak{f}$ is uniquely determined, up to a scalar coefficient, by $I=(\ell_1^{s+1},\ldots,\ell_{n}^{s+1}):\mathfrak{f}$ and the condition that no nonzero term of $\mathfrak{f}$ belongs to the ideal $(\ell_1^{s+1},\ldots,\ell_{n}^{s+1})$, one has $\mathfrak{f}=\lambda(\ell_1\cdots\ell_n)^{\mathfrak{g}}\mathfrak{f}_0$ for some nonzero $\lambda\in\kk.$ Hence, $(\ell_1\cdots\ell_n)^{\mathfrak{g}}$ divides $\mathfrak{f}.$  

Finally, we assert that $\mathfrak{g}$ is the largest index with this property, a claim that is obvious if $m_0=1,$ because in this case $\mathfrak{f}_0$ is a nonzero scalar. Thus, suppose $m_0\geq 2.$  If $\mathfrak{g}$ is not the largest index such that $(\ell_1\cdots\ell_n)^{\mathfrak{g}}$ divides $\mathfrak{f}$ then $\ell_1\cdots\ell_n$ divides $\mathfrak{f}_0.$ Hence, by Lemma~\ref{cond_iterado}, 
$$I=(\ell_1^{m_0},\ldots,\ell_n^{m_0}):\mathfrak{f}_0=(\ell_1^{m_0-1},\ldots,\ell_n^{m_0-1}):\frac{\mathfrak{f}_0}{\ell_1\cdots\ell_n},$$
so, $\{\ell_1^{m_0-1},\ldots,\ell_n^{m_0-1}\}\subset I,$  contradicting that $m_0$ is  least  such that $\{\ell_1^{m_0},\ldots,\ell_n^{m_0}\}\subset I.$
\qed

\section{The associated Rees algebra}\label{Rees_all}

\subsection{Equigenerated ideals of finite colength}

Let $I\subset R=\kk[x,y,z]$ be an equigenerated ideal of finite colength.
In this part  focus on the Rees algebra  $\mathcal{R}(I)\simeq R[It]\subset R[t]$, the associated graded ring ${\rm gr}_I(R)=\mathcal{R}(I)/I\mathcal{R}(I)$ and the fiber cone algebra $\mathcal{F}(I)=\mathcal{R}(I)/\fm \mathcal{R}(I)$, where $\fm:=(x,y,z)$.
The eventual goal is an application to the case where $I$ is Gorenstein.
The nature of the associated graded ring for Artinian Gorenstein rings in any dimension has been considered earlier by Iarrobino (\cite{Iarr}).

Note that the so-called condition $G_3$ is automatic since the ideal is $\fm$-primary. Some features in this section might have appeared elsewhere coming from a different angle.
Yet, it may be useful to have elementary proofs of the results below, where $G_3$ is not directly used.

\begin{Proposition}\label{main_colength} Let $R=\kk[x_1,\ldots,x_n]$ be a standard graded polynomial ring over a field and let $\fm$ be its maximal homogeneous ideal. Let $I$ be a $d$-equigenerated homogeneous $\fm$-primary ideal. Given an integer $m_0\geq 1$ such that $I^{m_0}=\fm^{dm_0}$, then the following hold:
	\begin{enumerate}
		\item[{\rm(a)}] $I^{m}=\fm^{md}$ for every $m\geq m_0.$
		\item[{\rm(b)}] The reduction number of $I$ is at most $\max\{m_0,r(\fm^d)\}$, where $r(\fm^d)$ denotes the reduction number of $\fm^d$.
		\item[{\rm(c)}] The {\rm (}regular{\rm )} rational map $\mathfrak{F}:\P^{n-1}\dasharrow\P^{\mu(I)-1}$ defined by  a set of forms spanning $I_d$ is birational onto the image.
		\item[{\rm(d)}] The Rees algebra $\mathcal{R}(I)$  satisfies the condition $R_{1}$ of Serre.
		\item[{\rm(e)}] $\depth {\rm gr}_I(R)=0.$
	\end{enumerate}
\end{Proposition}
\demo (a) One has $I^{m_0}\subset \fm^{d}I^{m_0-1}\subset \fm^{d}\fm^{(m_0-1)d}=I^{m_0}$, hence, $I^{m_0}=\fm^dI^{m_0-1}.$ Inducting on $m\geq m_0$, $$I^{m+1}=I^{m+1-m_0}I^{m_0}=I^{m+1-m_0}\fm^dI^{m_0-1}=I^{m}\fm^{d}=\fm^{(m+1)d}.$$

(b) Let $J\subset I$ be a homogeneous minimal reduction.
Since $\fm^d$ is the integral closure of $I$, then $J$ is also a minimal reduction of $\fm^d$. 
Setting $N=\max\{m_0,r(\fm^d)\},$ one has:
\begin{eqnarray}
	I^{N+1}&=&(\fm^{d})^{N+1}\quad \quad \quad\mbox{by (a)}\nonumber\\
	&=&J(\fm^{d})^{N}\quad\quad\quad\,\,\mbox{because $J$ is a minimal reduction of $\fm^d$}\nonumber\\
	&=&JI^{N} \quad \quad\quad\quad \,\,\,\,\mbox{by (a).}\nonumber
\end{eqnarray}

(c) By  (a),  the Hilbert polynomial $HP(\mathcal{F}(I),m)$ of the fiber cone $\mathcal{F}(I)$ is 

$$HP(\mathcal{F}(I),m)={md+n-1\choose n-1}=\frac{d^{n-1}}{(n-1)!}m^{n-1}+\mbox{lower degree terms of $m$}.$$ Hence, the multiplicity $e(\mathcal{F}(I))$ of  $\mathcal{F}(I)$  is $d^{n-1}.$ On the other hand, by \cite[Theorem 6.6 (a)]{Ram2} the degree $\deg(\mathfrak{F})$  of the rational map $\mathfrak{F}$ is
\begin{eqnarray}
	\deg(\mathfrak{F})&=& \frac{d^{n-1}}{e(\mathcal{F}(I))}. \nonumber
\end{eqnarray}
Thus, $\deg(\mathfrak{F})=1$, as asserted.

(d) Consider the Hilbert-Samuel polynomial ($m>\!\!> 0$)  $$\lambda(R/I^{m+1})=e_0(I){n+m\choose n}- e_1(I){n+m-1\choose n-1}+\mbox{lower degree terms of $m$}$$
and the Hilbert polynomial 
$$\lambda(R/\overline{I^{m+1}})=\overline{e}_0(I){n+m\choose n}- \overline{e}_1(I){n+m-1\choose n-1}+\mbox{lower degree terms of $m$}$$ 
where $\overline{I^{m+1}}$ denotes the integral closure of $I^{m+1}.$ By (a), $I^{m}=\overline{I^{m}}$ for every $m\geq m_0.$ Thus, in particular, $e_1(I)=\overline{e}_1(I).$ Hence, by \cite[Proposition 3.2]{syl2}, $\mathcal{R}(I)$ satisfies the condition $R_1$ of Serre.

(e) By (a), one has an exact sequence
$$0 \lar\mathcal{R}(I)  \lar \mathcal{R}(\fm^d) \lar C \lar 0,$$  
with  $C$  a module of finite length. In particular, $\depth C = 0.$ Since  $\mathcal{R}(\fm^d) $ is Cohen--Macaulay, then $\depth \mathcal{R}(I) = \depth C+1=1$. 

Now, clearly $\depth {\rm gr}_I(R)\leq \depth \mathcal{R}(I)=1$.
Supposing that $\depth {\rm gr}_I(R)>0$, let $a\in I\setminus I^2$ be such that its image in $I/I^2\subset {\rm gr}_I(R)$ is a regular element.
Then one has an exact sequence
$$ 0\rar {\rm gr}_I(R)(-1)\rar  \mathcal{R}(I)/a \mathcal{R}(I)\rar  \mathcal{R}_{R/(a)}(I/(a)) \rar 0$$
(see \cite[Proposition 5.1.11]{WolmBook1}).
Since $a$ is regular on $\mathcal{R}(I)$ then the middle term has depth zero, while the rightmost term -- being a Rees algebra over a Cohen--Macaulay ring of dimension $\geq 1$ -- has depth at least one.
It follows that ${\rm gr}_I(R)\simeq {\rm gr}_I(R)(-1)$ has depth zero; a contradiction. 
\qed

\subsection{Syzygetic ideals}

For the main result in this part recall the notion of a {\em syzygetic} ideal $I\subset R$ as being one such that the natural surjection $\mathcal{S}_R(I)\surjects \mathcal{R}_R(I)$ is an isomorphism in degree $\leq 2$. In particular, for such an ideal, $I^2$ coincides with the second symmetric power of $I$, hence the minimal number of generators of $I^2$ is given by ${{\mu(I)+1}\choose 2}$,  where $\mu(I)$ stands for the minimal number of generators of $I$.

In the ternary case we can bring over the fiber cone.

\begin{Theorem} \label{main_syzygetic}
	Let $I$ be a $d$-equigenerated $\fm$-primary homogeneous ideal in the standard polynomial ring $R=\kk[x,y,z]$, with $d\geq 2.$ If $I$ is syzygetic and minimally generated by $2d+1$ forms, the following hold:  
	\begin{enumerate}
		\item[{\rm(a)}] $I^2=\fm^{2d}.$ 
		\item[{\rm(b)}] $I^{m}=\fm^{md}$ for every $m\geq 2.$
		\item[{\rm(c)}] The reduction number of $I$ is $2.$
		\item[{\rm(d)}] The {\rm (}regular{\rm )} rational map $\mathfrak{F}:\P^{2}\dasharrow\P^{2d}$ defined by a set of forms spanning $I_d$ is birational onto the image.
		\item[{\rm(e)}] The Rees algebra $\mathcal{R}(I)$  satisfies the condition $R_{1}$ of Serre.
		\item[{\rm(f)}] $\depth {\rm gr}_I(R)=0.$
		\item[{\rm(g)}] The fiber cone $\mathcal{F}(I)$ is not Cohen-Macaulay.
	\end{enumerate}	
\end{Theorem}
\demo (a) Since $I$ is syzygetic then $$\mu(I^2)={\mu(I)+1\choose 2}={2d+2\choose 2}=\mu(\fm^{2d}).$$ Thus, $I^2\subset \fm^{2d}$ is a inclusion of homogeneous ideal generated in fixed degree $2d$ having the  same minimal number of homogenous generators. Hence, $I^2=\fm^{2d}.$

It remains now to prove items (c) and (g) because the others follow exactly as in Proposition~\ref{main_colength}.

(c) By (a) and Proposition~\ref{main_colength}, one has $r(I)\leq 2.$ On the other hand, since $I$ is syzygetic we have $2\leq r(I).$ Hence, $r(I)=2.$

(g) Suppose to the contrary. Then, by \cite[Proposition 1.2]{GST}, the reduction number $r(I)$ is the Castelnuovo-Mumford regularity ${\rm reg}(\mathcal{F}(I))$ of $\mathcal{F}(I).$
By (c), the latter is $2$.
But since $I$ is syzygetic, the defining ideal of $\mathcal{F}(I)$ over $S:=\kk[T_1,\ldots,T_{2d+1}]$ admits no forms of degree $2$, hence is generated in the single degree $3$. In particular, the minimal graded free resolution of $\mathcal{F}(I)$ over $S $ is linear:
$$0\to S(-N+1)^{\beta_{N-3}}\to\cdots\to S(-3)^{\beta_1}\to S.$$
By \cite[Theorem 1.2]{HM}, the multiplicity of the fiber cone $\mathcal{F}(I)$ is
$$e(\mathcal{F}(I))={\mu(I)-1\choose2}={2d \choose 2}.$$
Now consider the rational map $\mathfrak{F}:\P^2\dasharrow \P^{2d}$ defined by the given generators of $I$ in degree $d$, and let $\deg (\mathfrak{F})$ denote the degree of $\mathfrak{F}$.
Since $I$ is equigenerated then $\mathcal{F}(I)$ is isomorphic to the $\kk$-subalgebra of $R$ generated by the vector space $I_d$, while the latter is up to degree normalization the homogeneous defining ideal of the image of $\mathfrak{F}$.
Then, by \cite[Theorem 6.6 (a)]{Ram2} one has 
$${2d \choose 2}\deg (\mathfrak{F})=e(\mathcal{F}(I)) \deg (\mathfrak{F})=d^2.$$
Since $\mathfrak{F}$ is birational,
$2d-1=d$, which is absurd for $d\geq 2$.
\qed

\subsection{Application to the Gorenstein case}

\begin{Corollary}\label{fiber_cone_of_Gorenstein} {\rm (char$(\kk)\neq 2$)}
	Let $I$ denote a codimension $3$ homogeneous Gorenstein ideal in $\kk[x,y,z]$ with datum $(d,2d+1)$, where $d\geq 2$. Then all assertions of {\rm Theorem~\ref{main_syzygetic}} hold true.
\end{Corollary}
\demo Since char$(\kk)\neq 2$, then $I$ is syzygetic (\cite[Proposition 2.8]{Trento}). 
\qed 

\smallskip

For the non-linear case, we have the following:

\begin{Proposition}
	Let $I\subset R=\kk[x,y,z]$ be a codimension $3$ Gorenstein ideal with datum $(d,r)$ and skew degree $d'.$ Let $\mathfrak{F}:\mathbb{P}^2\dasharrow\mathbb{P}^{r-1}$ be the rational map defined by the linear system $I_d.$ If the reduction number of $I$ is at most $2$ and  $\mathcal{F}(I)$ is Cohen-Macaulay then:
	\begin{enumerate}
		\item[{\rm (a)}]  $(r-2)$ divides $d'^2.$
		\item[{\rm(b)}] If $r\geq 5$  then $\mathfrak{F}$  is not birational onto the image.
	\end{enumerate}
\end{Proposition}
\demo (a) Since $I$ is syzygetic, the assumption implies that the reduction number of $I$ is exactly $2$.
Since $\mathcal{F}(I)$ is Cohen-Macaulay then the same argument as in the proof of Theorem~\ref{main_syzygetic} (g) yields  $e(\mathcal{F}(I))={r-1\choose 2}.$  Again, by \cite[Theorem 6.6]{Ram2}, ${r-1\choose 2} \deg (\mathfrak{F})=d^2.$ By definition, $d=(r-1)d'/2$,  hence
\begin{equation}\label{degF}
	2(r-2)\deg (\mathfrak{F})=(r-1)d'^{2}.
\end{equation}
Since gcd$\{(r-2),(r-1)\}=1$ then $(r-2)$ divides $d'^2$, as desired.

(b) Since $(r-1)/2> 1$ then \eqref{degF}forces $\deg (\mathfrak{F})>1.$ Hence, $\mathfrak{F}$ is not birational.
\qed

\section{Appendix}

In this part we give some explicit results, which can be used in particular to give an alternative proof of Theorem~\ref{Old_Theorem 2.7} without using either \cite{HoLa} or \cite{Stan}.

Assume the following setup:

\begin{Setup}\label{setup} 
	
	Set  $R=\kk[x,y,z]$.
	
	$\bullet$ $r\geq 5 $ is a odd integer and $d:=(r-1)d'/2$ for some integer $d'\geq 1.$
	
	$\bullet$ $e:=d+d'-3$.
\end{Setup}

\begin{Proposition}\label{examplemonomialbis} With the data of {\rm Setup~\ref{setup}},  let $I_{d',r}\subset R=\kk[x,y,z]$ be   the ideal $$I_{d',r}=(z^{d'}(x^{d'},y^{d'})^{(r-3)/2},x^{2d'}(x^{d'},y^{d'})^{(r-5)/2},y^d,z^d).$$ 
	Then,
	\begin{enumerate}
		\item[{\rm(a)}] The minimal graded free resolution of $R/I_{d',r}$ has the form
		{\small$$0\to\begin{array}{c}R(-(d+2d'))^{(r-5)/2}\\\oplus\\R(-(d+3d'))\\\oplus\\R(-(2d))^{(r-3)/2}\end{array}\to\begin{array}{c}R(-(d+d'))^{\frac{3(r-1)-6}{2}}\\\oplus\\ R(-(d+2d'))\\\oplus\\ R(-(2d-d'))^{(r-1)/2}\end{array}\to R(-d)^{r}\to R\to R/I_{d',r}\to 0$$}
		\item[{\rm(b)}] The decomposition structure of the socle of $I_{d',r}$ is 
		{\footnotesize \begin{equation*}
				{\rm Soc}(R/I_{d',r})=\left\{\begin{array}{cc}
					\kk(-(4d'-3))\oplus\kk(-(5d'-3)),&\mbox{if}\,\, r=5\\
					\kk(-(d+2d'-3))^{(r-5)/2}\oplus\kk(-(d+3d'-3))\oplus\kk(-(2d-3))^{(r-3)/2},&\mbox{if}\,\,r\geq 7
				\end{array}\right.
		\end{equation*}}
		\item[{\rm(c)}] If, moreover, $r\geq 7$ then  $V(I_{D}(M_{d,r,d+d'-3}))$  is a proper subset of $(\P^N)^r$, where $D:=\dim_{\kk}(R_{2d+d'-3})$.	
	\end{enumerate}
\end{Proposition}
\demo
 (a) We will first prove the case $d'=1$ (in particular, $r=2d+1$).  For convenience, we  write the generators of $I$ in the following block shaped matrix:

\vskip-10pt

\begin{equation}\label{phi_0}
\phi_0:=\left[\begin{array}{cccc}z\ff&x^2\gg& y^{d} &z^d\end{array}\right],
\end{equation}
where $\ff:=[x^{d-1}\; x^{d-2}y\;\cdots\;y^{d-1}]$ and $\gg:=[x^{d-2}\;x^{d-3}y\;\cdots\;y^{d-2}]$.
Then $\phi_0$ is the matrix of the map $R(-d)^{2d+1}\to R$.

Throughout  $ \ee_j $ and $\boldsymbol0_{i \times j} $ will denote the identity matrix of order $j$  and the $i\times j$ null matrix, respectively.
Consider further the following matrices, which will be candidates to first and second syzygies.

\begin{equation}\label{phi_1} \phi_1=\left[\begin{array}{ccccc|c|c}
\aa_1&\aa_2&\aa_3&\aa_4&\boldsymbol0_{d\times (d-2)}&\boldsymbol{0}_{d\times 1}&\aa_7\\
\bb_1&\boldsymbol0_{(d-1)\times1}&\bb_3&\boldsymbol0_{(d-1)\times 1}&\bb_5&\bb_6&\boldsymbol0_{(d-1)\times d}\\
\boldsymbol{0}_{1\times(d-2)}&\boldsymbol{0}_{1\times1}&\boldsymbol0_{1\times(d-1)}&\cc_4&\boldsymbol0_{1\times(d-2)}&\cc_6&\boldsymbol{0}_{1\times d}\\
\boldsymbol0_{1\times(d-2)}&\boldsymbol0_{1\times1}&\boldsymbol0_{1\times (d-1)}&\boldsymbol0_{1\times 1 }&\boldsymbol0_{1\times (d-2)}&0_{1\times 1}&\ddd_7
\end{array}\right],
\end{equation}
with three blocks in standard degrees $1, 2$ and $d-1$ from left to right, respectively -- if needed we will refer to these three blocks as {\em degree blocks}.

Here
{\small$$\aa_1=\left[\begin{array}{c}y\ee_{d-2}\\\boldsymbol0_{2\times (d-2)}\end{array}\right],\,\, \aa_{2}=\left[\begin{array}{c} \boldsymbol{0}_{(d-2)\times 1}\\-y\\x\end{array}\right],\,\, \aa_3=\left[\begin{array}{c}x\ee_{d-1}\\\boldsymbol0_{1\times (d-1)}\end{array}\right],\,\, \aa_{4}=\left[\begin{array}{c}\boldsymbol0_{(d-1)\times 1}\\y\end{array}\right],\,\aa_{7}=-z^{d-1}\ee_d,$$}

{\small$$\quad\bb_{1}=\left[\begin{array}{c}\boldsymbol0_{1\times(d-2)}\\ -z\ee_{d-2}\end{array}\right],\quad\bb_3=-z\ee_{d-1},\quad\bb_{5}=\left[\begin{array}{ccccc}
	-y&0&\cdots&0\\
	x&-y&\cdots&0\\
	\vdots&\ddots&\ddots&\vdots\\
	0&\cdots&x&-y\\
	0&\cdots&0&x\\
	\end{array}\right],\quad \bb_{6}=\left[\begin{array}{c}\boldsymbol{0}_{(d-2)\times 1}\\\hline-y^2\end{array}\right],$$}

{$$ \cc_4=\left[\begin{array}{c}-z\end{array}\right],\quad\cc_6=\left[\begin{array}{c}x^2\end{array}\right], \quad\ddd_7=\left[\begin{array}{cccccc}
	x^{d-1}&x^{d-2}y&\cdots&y^{d-1}
	\end{array}\right].$$}

The second matrix is
\begin{equation}\label{phi_2}
\,\phi_{2}=\left[\begin{array}{cccc}
\ss_1&\boldsymbol0_{(d-2)\times 1}&\uu_1&\boldsymbol0_{(d-2)\times 1}\\
\boldsymbol0_{1\times(d-2)}&\tt_{2}&\boldsymbol0_{1\times(d-2)}&\vv_2\\
\ss_{3}&\tt_3&\uu_3&\boldsymbol0_{(d-1)\times 1}\\
\boldsymbol0_{1\times d-2}&\tt_4&\boldsymbol0_{1\times(d-2)}&\boldsymbol0_{1\times 1}\\
\ss_5&\boldsymbol0_{(d-2)\times 1}&\boldsymbol0_{(d-2)\times(d-2)}&\boldsymbol0_{(d-2)\times1}\\
\boldsymbol0_{1\times (d-2)}&\tt_6&\boldsymbol0_{1\times(d-2) }&\boldsymbol0_{1\times 1}\\
\boldsymbol0_{d\times(d-2)}&\boldsymbol0_{d\times 1}&\uu_{7}&\vv_7\end{array}\right],
\end{equation}
where

{$$\ss_1=x\ee_{d-2},\quad\ss_3=\left[\begin{array}{c}-y\ee_{d-2}\\\boldsymbol0_{1\times(d-2)}\end{array}\right],\quad\ss_5=z\ee_{d-2},$$}
{$$\tt_2=\left[\begin{array}{c}-xy\end{array}\right],\quad\tt_3=\left[\begin{array}{c}\boldsymbol0_{(d-2)\times 1}\\-y^2\end{array}\right],\quad\tt_4=\left[\begin{array}{c}x^2\end{array}\right],\quad\tt_6=\left[\begin{array}{c}z\end{array}\right],$$}
{$$\uu_1=-z^{d-1}\ee_{d-2},\quad\uu_3=\left[\begin{array}{c}\boldsymbol{0}_{1\times(d-2)}\\z^{d-1}\ee_{d-2}\end{array}\right],\quad\uu_7=\left[\begin{array}{ccccc}
	-y&0&\cdots&0\\
	x&-y&\cdots&0\\
	\vdots&\ddots&\ddots&\vdots\\
	0&\cdots&x&-y\\
	0&\cdots&0&x\\
	0&\ldots&0&0
	\end{array}\right],$$}
{$$\vv_2=\left[\begin{array}{c}z^{d-1}\end{array}\right],\quad \vv_7=\left[\begin{array}{c}\boldsymbol{0}_{(d-2)\times 1}\\-y\\x\end{array}\right].$$}

The goal of this seemingly bizarre block wise way of writing matrices is to adjust page fitting and reading easiness.

\smallskip

\noindent{\sc Claim 1.} The sequence of $R$-maps
\begin{equation}\label{isacomplex}
0\to R^{2d-2}\stackrel{\phi_2}\longrightarrow R^{4d-2}\stackrel{\phi_1}\longrightarrow R^{2d+1}\stackrel{\phi_0}\longrightarrow R
\end{equation}
is a complex of $R$-modules.

The fact that the composite 
{\small$$\phi_0\cdot \phi_1=\left[\begin{array}{ccccccc}
	z\ff\aa_1+x^2\gg\bb_1&z\ff\aa_2&z\ff\aa_3+x^2\gg\bb_3&z\ff\aa_4+y^d\cc_4&x^2\gg\bb_5&x^2\gg\bb_6+y^d\cc_6&z\ff\aa_7+z^d\ddd_7\end{array}\right]$$}
is a null matrix is a routine exercise in syzygies of monomials as reduced Koszul relations. Yet, the shape of $\phi_1$ will be of relevance later for rank and minors computation.

The other composite

{\small $$\phi_1\cdot\phi_2=\left[\begin{array}{ccccccc}
	\aa_1\ss_1+\aa_3\ss_3&\aa_2\tt_2+\aa_3\tt_3+\aa_4\tt_4&\aa_1\uu_1+\aa_3\uu_3+\aa_7\uu_7&\aa_2\vv_2+\aa_7\vv_7\\
	\bb_1\ss_1+\bb_3\ss_3+\bb_5\ss_5&\bb_3\tt_3+\bb_6\tt_6&\bb_1\uu_1+\bb_3\uu_3&\boldsymbol0_{(d-1)\times 1}\\
	\boldsymbol0_{1\times (d-2)}&\cc_4\tt_4+\cc_6\tt_6&\boldsymbol0_{1\times(d-2)}&\boldsymbol0_{1\times1}\\
	\boldsymbol0_{1\times(d-2)}&\boldsymbol0_{1\times 1}&\ddd_7\uu_7&\ddd_7\vv_7
	\end{array}\right]$$
}
is a bit more delicate, but all calculations are straightforward.

\medskip

\noindent{\sc Claim 2.} The complex \eqref{isacomplex} is acyclic.

For the argument we use the Buchsbaum--Eisenbud acyclicity criterion. Obviously, $\rk\phi_0= 1$ and  $\hht I_1(\phi_0)\geq 1.$

We next focus on $\phi_1$, aiming at showing that $\rk\phi_1\geq 2d$ (hence, $\rk\phi_1= 2d$) and $\hht I_{2d}(\phi_1)\geq 2.$

For this, we single out the following $2d\times 2d$ submatrices.

$$A:=\left[\begin{array}{cccc|c}
\aa_3&\aa_4&\boldsymbol0_{d\times(d-2)}&\boldsymbol{0}_{d\times 1}&\\ 
\bb_3&\boldsymbol0_{(d-1)\times 1}&\bb_5&\bb_6&\boldsymbol\rho\\ 
\boldsymbol0_{1\times(d-1)}&\boldsymbol0_{1\times1}&\boldsymbol0_{1\times(d-2)}&\boldsymbol0_{1\times1}& 
\end{array}\right],$$
formed with rows $1,2,\ldots,2d-1,2d+1$ and columns $d,d+1, \ldots, 3d-1$, where 
$$\boldsymbol\rho=(-z^{d-1}\; 0\;\ldots \; 0\; x^{d-1})^t$$ 
is the first column of the rightmost degree block of $\phi_1$,
and

$$B=\left[\begin{array}{cccccccc}
\aa_3&\aa_{4}&\aa_7\\
\bb_3&\boldsymbol{0}_{(d-1)\times 1}&\boldsymbol{0}_{(d-1)\times d}\\
\boldsymbol{0}_{1\times(d-1)}&\cc_4&\boldsymbol{0}_{1\times d}
\end{array}\right],$$
formed with rows $1,2,\ldots,2d$ and columns $d,d+1,\ldots,2d,3d,\ldots, 4d-2$.

Expanding conveniently, we have
{\small\begin{eqnarray}\label{detA}\det A&=&x^{d-1}\det\left[\begin{array}{ccccc}
	\aa_3&\aa_4&\boldsymbol0_{d\times(d-2)}&\boldsymbol{0}_{d\times 1}\\
	\bb_3&\boldsymbol0_{(d-1)\times 1}&\bb_5&\bb_6
	\end{array}\right]\nonumber\\
	&=&x^{d-1}\det\left[\begin{array}{ccc}
	\aa_3&\aa_4
	\end{array}\right]\det\left[\begin{array}{ccc}
	\bb_5&\bb_6
	\end{array}\right]\nonumber\\
	&=&x^{d-1}\det\left[\begin{array}{cc}x\ee_{d-1}&\boldsymbol0_{(d-1)\times 1}\\\boldsymbol0_{1\times (d-1)}&y\end{array}\right]\det\left[\begin{array}{ccccc}
	-y&0&\cdots&0&0\\
	x&-y&\cdots&0&0\\
	\vdots&\ddots&\ddots&\vdots&\vdots\\
	0&\cdots&x&-y&0\\
	0&\cdots&0&x&-y^2\\
	\end{array}\right]\nonumber\\
	&=&x^{d-1}(yx^{d-1})((-1)^{d-1}y^d)=(-1)^{d-1}x^{2d-2}y^{d+1}\in I_{2d}(\phi_1)
	\end{eqnarray}}
and

{\small
	\begin{eqnarray}\label{detB}
	\det B&=&(-1)^d \cc_4\det\aa_7\det\bb_3=(-1)^d (-z) (-1)^{d}z^{(d(d-1))} (-1)^{d-1}(z^{d-1})\\
	&=&(-1)^{3d} z^{d^2}\in I_{2d}(\phi_1).\nonumber
	\end{eqnarray}}
Therefore, we are  through.

\smallskip

Next get to $\phi_2$, for which we want to prove that $\rk\phi_2\geq 2d-2$ and $\hht I_{2d-2}(\phi_2)\geq 3$.
Note that, since $\rk\phi_1=2d$, then $\rk\phi_2\leq 2d-2$, so we have derived the sought equality.

The determinants of the following three $(2d-2)\times(2d-2)$ submatrices of $\phi_2$ will be shown to form a regular sequence:

{$$S=\left[\begin{array}{ccccc}\ss_1&\boldsymbol0_{(d-2)\times 1}&\uu_1&\boldsymbol0_{(d-2)\times 1}\\
	\boldsymbol0_{1\times (d-2)}&\tt_4&\boldsymbol0_{1\times(d-2)}&\boldsymbol0_{1\times 1}\\
	\boldsymbol0_{(d-1)\times(d-2)}&\boldsymbol0_{(d-1)\times1}&\tilde{\uu}_7&\tilde{\vv}_7\end{array}\right],$$
	
	$$T=\left[\begin{array}{cccccc}\ss_3&\tt_3&\uu_3&\boldsymbol0_{(d-1)\times1}\\\boldsymbol0_{(d-1)\times(d-2)}&\boldsymbol0_{(d-1)\times 1}&\bar{\uu}_7&\bar{\vv}_{7}\end{array}\right]$$
	and
	
	$$U=\left[\begin{array}{cccccc}
	\ss_1&\boldsymbol0_{(d-2)\times 1}&\uu_1&\boldsymbol0_{(d-2)\times1}\\
	\boldsymbol0_{1\times(d-2)}&\tt_2&\boldsymbol0_{1\times(d-2)}&\vv_2\\
	\ss_5&\boldsymbol0_{(d-2)\times1}&\boldsymbol0_{(d-2)\times(d-2)}&\boldsymbol0_{(d-2)\times 1}\\
	\boldsymbol0_{1\times(d-2)}&\tt_6&\boldsymbol0_{1\times(d-2)}&\boldsymbol0_{1\times1}\end{array}\right].$$
	Here
	\begin{enumerate}
		\item[$\bullet$] $\tilde{\uu}_7$ and $\tilde{\vv}_7$ are the submatrices obtained from $\uu_7$ and $\vv_7,$ respectively, by omitting the first row.
		\item[$\bullet$] $\overline{\uu}_7$ and $\overline{\vv}_7$ are the submatrices obtained from $\uu_7$ and $\vv_7,$ respectively, by omitting the last row.
	\end{enumerate}
	
	The calculation is straightforward:
	
	{\small\begin{eqnarray}
		\label{detS}\det S&=&\tt_4\det\ss_1\det\left[\begin{array}{cc}
		\tilde{\uu}_7&\tilde{\vv}_7
		\end{array}\right]\nonumber\\
		&=& \tt_4\det\ss_1\det \left[\begin{array}{ccccc}
		x&-y&\cdots&0&0\\
		\vdots&\ddots&\ddots&\vdots&\vdots\\
		0&\cdots&x&-y&0\\
		0&\cdots&0&x&-y\\
		0&\ldots&0&0&x
		\end{array}\right]\nonumber\\
		&=&x^2\cdot x^{d-2}\cdot x^{d-1}=x^{2d-1},
		\end{eqnarray}}
	
	{\small\begin{eqnarray*}
		\det T&=&\det \left[\begin{array}{cc}\ss_3&\tt_3\end{array}\right]\det \left[\begin{array}{cc}\bar{\uu}_7&\bar{\vv}_7\end{array}\right]\nonumber\\
		&=&\det \left[\begin{array}{cc}\begin{array}{c}-y\ee_{d-2}\\\boldsymbol0_{1\times(d-2)}\end{array}&\begin{array}{c}\boldsymbol0_{(d-2)\times 1}\\-y^2\end{array}\end{array}\right]\det \left[\begin{array}{ccccc}
		-y&0&\cdots&0&0\\
		x&-y&\cdots&0&0\\
		\vdots&\ddots&\ddots&\vdots&\vdots\\
		0&\cdots&x&-y&0\\
		0&\cdots&0&x&-y
		\end{array}\right]\nonumber\\
		&=&y^{2d-1},\end{eqnarray*}}
	{\small \begin{eqnarray*}
		\det U&=&\pm\det \ss_5\cdot\det\left[\begin{array}{cccccc}
		\boldsymbol0_{(d-2)\times 1}&\uu_1&\boldsymbol0_{(d-2)\times1}\\
		\tt_2&\boldsymbol0_{1\times(d-2)}&\vv_2\\
		\tt_6&\boldsymbol0_{1\times(d-2)}&\boldsymbol0_{1\times1}\end{array}\right]\nonumber\\
		&=&\pm\det\ss_5\det\uu_1\det\vv_2\det\tt_6=\pm z^{d-2}(-1)^dz^{(d-1)(d-2)}z^{d-1}z=\pm (-1)^dz^{d(d-1)}.
		\end{eqnarray*}}
	This completes the argument on the acyclicity criterion, hence also the proof for the case $d'=1.$
	
	Now, deal with the case of arbitrary $d'>1.$ Consider the endomorphism of $\kk$-algebras $$\zeta:\kk[x,y,z]\to \kk[x,y,z],\quad x\mapsto x^{d'},\,y\mapsto y^{d'},\,z\mapsto z^{d'}.$$ Note that $I_{d',r}$ is the extension of $I_{1,r}$ by the endomorphim $\zeta$. Consider the matrices $\phi_0,$ $\phi_1$ and $\phi_2$ as in \eqref{phi_0}, \eqref{phi_1} and \eqref{phi_2}. Let $\zeta(\phi_i)$ ($i=0,1,2$) denote the matrix obtained from $\phi_i$ by evaluating $\zeta$ in each of its entries. Obviously,
	$$\zeta(\phi_0)\cdot \zeta(\phi_1)=0\quad\mbox{and}\quad \zeta(\phi_1)\cdot \zeta(\phi_2)=0.$$
	Hence, the sequence
	\begin{equation}\label{isacomplexbis}
	0\to R^{r-3}\stackrel{\zeta(\phi_2)}\longrightarrow R^{2r-4}\stackrel{\zeta(\phi_1)}\longrightarrow R^{r}\stackrel{\zeta(\phi_0)}\longrightarrow R
	\end{equation}
	is a complex.
	
	As shown above, there are $(r-1)\times (r-1)$ submatrices $A$ and $B$ of $\phi_1$  such that
	$$\det A= (-1)^{d-1}  x^{r-3} y^{(r+1)/2} \in I_{2d}(\phi_1)\quad\mbox{and}\quad \det B=  (-1) ^{3d} z^{((r-1)/2)^2}\in I_{2d}(\phi_1).$$
	Thus, 
	$$\det \zeta(A)=(-1)^{d-1}  x^{(r-3)d'} y^{((r+1)/2)d'}\in I_{r-1}(\zeta(\phi_1)),\, \det\zeta(B)=(-1) ^{3d} z^{((r-1)/2)^2d'}\in I_{r-1}(\zeta(\phi_1)).$$
	Hence,
$	\rk\zeta(\phi_1)= 2d$ and $\hht I_{r-1}(\zeta(\phi_1))\geq 2.$

	The proof for $d'=1$ also guarantees the existence of $(r-3)\times(r-3)$ submatrices $S,T,U$ of $\phi_2$  such that
	$$\det S=x^{r-2},\,\,\det T= y^{r-2},\, \,\det U=\pm (-1)^d z^{(r-1)(r-3)/4}.$$
	Thus, $\det \zeta(S)= x^{(r-2)d'},\,\det \zeta(T)= y^{(r-2)d'},\, \det \zeta(U)=\pm (-1)^d  z^{((r-1)(r-3)/4)d'},$ hence
	\begin{equation}\label{criterionzetaphi2}
	\rk\zeta(\phi_2)= r-3\quad \mbox{and}\quad \hht I_{r-3}(\zeta(\phi_2))\geq 3.
	\end{equation}
	Therefore, the complex \ref{isacomplexbis} is acyclic. 
	
	By construction, the cokernel of $\zeta(\phi_0)$ is $R/I_{d',r}.$ So,
	
	\begin{equation}\label{resRI'}
	0\to R^{r-3}\stackrel{\zeta(\phi_2)}\longrightarrow R^{2r-4}\stackrel{\zeta(\phi_1)}\longrightarrow R^{r}\stackrel{\zeta(\phi_0)}\longrightarrow R\to R/I_{d',r}\to 0
	\end{equation} 
	is a free resolution of $R/I_{d',r}.$ Finally, by observing the degrees of the entries of the matrices $\zeta(\phi_0),$  $\zeta(\phi_1)$ and $\zeta(\phi_2),$ the resolution \eqref{resRI'} we see that it is a minimal graded free resolution for $R/I_{d',r}$ as stated in the proposition.
	
	\smallskip
	
	(b) This is a consequence of  (a) via a well-known argument (see, e.g., \cite[Lemma 1.3]{KU}).
	
	\smallskip
	
	(c) Let $f_1,\ldots,f_{r}$ denote the given set of generators of $I_{d',r}.$ Then:
	\begin{eqnarray}
	R_{2d+d'-3}&=&I_{2d+d'-3}\quad(\mbox{by item (b)} )\nonumber\\ 
	&=&R_{d+d'-3}I_{d}\nonumber\\
	&=&R_{e}f_1+\cdots+R_{e}f_{r}.
	\end{eqnarray}
	Hence, by property \ref{B}, $P_{\ff}\notin V(I_{D}(M_{d,r,e})).$ In particular, $V(I_{D}(M_{d,r,e}))$
	is a proper  subset of $(\P^N)^r.$
	\qed
	
	\bigskip

In order to overcome the obstruction $r\geq 7$ in Proposition~\ref{examplemonomialbis} (c), we introduce the  following particular construct in five generators.

\begin{Proposition}\label{examplemonomialbisbis} Given an integer $d'\geq 2$, consider the following ideal of $R=\kk[x,y,z]$
	$$I=(x^{2d'},y^{2d'},z^{2d'},x^{d'}y^{d'},xz^{2d'-1}).$$
	Then:
	\item[{\rm(a)}] The minimal graded free resolution of $R/I$ has the form:
	{\small\begin{equation}
			0\to\begin{array}{c}R(-(4d'+1))\\\oplus\\R(-(5d'-1))^2\end{array}\to\begin{array}{c}R(-(2d'+1))\\\oplus\\R(-3d')^2\\\oplus\\R(-(4d'-1))^2\\\oplus\\R(-4d')^2\end{array}\to R(-2d')^5\to R\to  R/I\to0
	\end{equation}}
	\item[{\rm(b)}] The structure decomposition of the socle of $R/I$  is 
	\begin{equation}
		{\rm Soc}(R/I)=\kk(-(4d'-2))\oplus \kk(-(5d'-4))^2
	\end{equation}
	\item[{\rm(c)}] $V(I_{D}(M_{2d',5,3d'-3}))$  is a proper  subset of $(\P^N)^5$, where $D:=\dim_{\kk}(R_{d+3d'-3})$.	
\end{Proposition} 
\demo
Once more, write $\phi_0=[x^{2d'}\,y^{2d'}\,z^{2d'}\,x^{d'}y^{d'}\,xz^{2d'-1}]$
	and introduce first and second syzygies candidates:

	$$\phi_1=\left[\begin{array}{ccccccc}
	0&0&-y^{d'}&0&-z^{2d'-1}&0&0\\
	0&x^{d'}&0&0&0&-z^{2d'}&-xz^{2d'-1}\\
	x&0&0&0&0&y^{2d'}&0\\
	0&-y^{d'}&x^{d'}&-z^{2d'-1}&0&0&0\\
	-z&0&0&x^{d'-1}y^{d'}&x^{2d'-1}&0&y^{2d'}
	\end{array}\right]$$
	{and}
	$$\phi_2=\left[\begin{array}{ccc}y^{2d'}&0&0\\0&0&z^{2d'-1}\\0&z^{2d'-1}&0\\0&x^{d'}&-y^{d'}\\0&-y^{d'}&0\\-x&0&0\\z&0&x^{d'-1}\end{array}\right].$$
	A straightforward calculation yields
	$$\phi_0\cdot\phi_1=0\quad\mbox{and}\quad\phi_1\cdot\phi_2=0.$$
	Thus, the sequence $$0\to R^{3}\stackrel{\phi_2}\to R^{7}\stackrel{\phi_1}\to R^5\stackrel{\phi_0}\to R$$ is a complex.

	Now consider the following $4\times4$ submatrices of $\phi_1$
	{$$A=\left[\begin{array}{ccccc}
		0&0&-z^{2d'-1}&0\\
		0&0&0&-z^{2d'}\\
		0&-z^{2d'-1}&0&0\\
		-z&x^{d'-1}y^{d'}&x^{2d'-1}&0
		\end{array}\right]\quad\mbox{and}\quad 
		B=\left[\begin{array}{cccc}
		0&x^{d'}&0&0\\
		x&0&0&0\\
		0&-y^{d'}&x^{d'}&0\\
		-z&0&0&x^{2d'-1}
		\end{array}\right].$$}
	As easily seen, $\det A= -z^{6d'-1}$ and $\det B=-x^{4d'}.$ Therefore, $\rank\phi_1=4$ and $\hht I_4(\phi_1)\geq 2$.
	
	On the other hand, for the following $3\times 3$ submatrices of $\phi_2$ 
	
	{\small$$S=\left[\begin{array}{ccc}
		y^{2d'}&0&0\\
		0&x^{d'}&-y^{d'}\\
		0&-y^{d'}&0
		\end{array}\right],\quad
		T=\left[\begin{array}{ccc}
		0&0&z^{2d'-1}\\
		0&z^{2d'-1}&0\\
		z&0&x^{d'-1}
		\end{array}\right]\quad\mbox{and}\quad
		U=\left[\begin{array}{ccc}0&x^{d'}&-y^{d'}\\-x&0&0\\z&0&x^{d'-1}\end{array}\right]
		$$}
	one has $\det S=-y^{4d'},$ $\det T=-z^{4d'-1}$ and $\det U=x^{2d'},$ forming a regular sequence. 
	Thus, (a) holds.
	
	(b) and (c) follow easily as before.
	\qed
	
	\medskip
	
	We now give an alternative proof of Theorem~\ref{Old_Theorem 2.7}.
	
	\begin{Theorem}{\rm (Theorem~\ref{Old_Theorem 2.7} (bis))}\label{second_proof}
		Let  $I\subset R=\kk[x,y,z]$ be a Gorenstein ideal generated by a general set of $r\geq 5$  forms of degree $d\geq 2$. Then $r=5$ and $d=2$.
	\end{Theorem}
	\demo
	Since $I$ is generated by a general set of sufficiently many forms, it has finite colength (see Lemma~\ref{extremal_properties} (iv)).
	Let $\{f_1,\ldots,f_r\}$ be a set of such forms of degree $d$.  
	Since the socle degree of $I$ is $2d+d'-3$ (see Lemma~\ref{socle_degree_general}), then $R_{d+d'-3}f_1+\cdots+R_{d+d'-3}f_r\subsetneq R_{2d+d'-3}.$
	Therefore, \ref{B} implies that $P_{\ff}\in V(I_{D}(M_{d,r,d+d'-3}))=V(I_{D}(M_{((r-1)/2)d',r,((r-1)/2)d'+d'-3}))$, where $D=\dim_{\kk}(R_{2d+d'-3})$.
	
	Suppose that  $d\geq 3$. Then $((r-1)/2)d'=d\geq 3$, i.e., $(r-1)d'\geq 6$. Now, either $d'=1$ and $r\geq 7$, or else $d'\geq 2$.  Thus, in any case,  either Proposition~\ref{examplemonomialbis} or Proposition~\ref{examplemonomialbisbis} implies that $V(I_{D}(M_{d,r,d+d'-3}))$ is a proper Zariski closed subset of $(\P^N)^r.$  
	This shows that all Gorenstein ideals of codimension $3$, of given degree $d\geq 3$ and $r\geq 5$ number of generators, are parameterized by a proper Zariski closed subset of the space of parameters.  
	\qed


{\footnotesize
\begin{center} {\bf Addresses:} \end{center}

\medskip

\begin{tabular}{ll}
{\sc Dayane  Lira} & \hspace{3cm} {\sc Zaqueu Ramos}\\
Departamento de Matem\'atica, CCEN & \hspace{3cm} Departamento de Matem\'atica, CCET\\
Universidade Federal da Para\'iba & \hspace{3cm} Universidade Federal de Sergipe\\
 58051-900 J. Pessoa, PB, Brazil & \hspace{3cm} 49100-000 S\~ao Cristov\~ao, Sergipe, Brazil\\
 dayannematematica@gmail.com & \hspace{3cm} zaqueu@mat.ufs.br
\end{tabular}

\smallskip

\begin{center}
\noindent {\sc Aron Simis}\\
Departamento de Matem\'atica, CCEN\\ 
Universidade Federal de Pernambuco\\ 
50740-560 Recife, PE, Brazil\\
{\em e-mail}:  aron@dmat.ufpe.br
\end{center}
}

\end{document}